\documentclass[journal,twoside]{IEEEtran}
\usepackage{amsmath,amsfonts}
\usepackage{algorithmic}
\usepackage{algorithm}
\usepackage{array}
\usepackage[caption=false,font=normalsize,labelfont=sf,textfont=sf]{subfig}
\usepackage{textcomp}
\usepackage{stfloats}
\usepackage{url}
\usepackage{verbatim}
\usepackage{amssymb}
\usepackage{mathtools}
\usepackage{graphicx}
\usepackage{cite}
\usepackage{hyperref}
\usepackage{xcolor}
\usepackage{multirow}
\usepackage{todonotes}
\usepackage{amsthm}
\newtheorem{lemma}{Lemma}

\newcommand{\X}{\mathbb{X}}
\newcommand{\Y}{\mathbb{Y}}
\newcommand{\R}{\mathbb{R}} 
\DeclareMathOperator*{\argmin}{arg\,min}

\newtheorem{algo}[lemma]{Algorithm}
\newcommand{\genPar}{\zeta}
\newcommand{\regPar}{\chi}
\newcommand{\dt}{\,\mathrm{d}t}
\newcommand{\dw}{\,\mathrm{d}w}
\newcommand{\dynR}{\mathcal{R}}

\begin{document}

\title{Transformer Causality Regularization for \\ Dynamic Inverse Problems}

\author{Gesa Sarnighausen, Anne Wald, Andreas Hauptmann~\IEEEmembership{Senior Member,~IEEE}
\thanks{AH is supported in part by the Research council of
Finland (Flagship
of Advanced Mathematics for Sensing Imaging and Modelling (FAME), Project No. 359186; Academy Research Fellow, Project No. 338408, and Academy Project No. 370528). GS and AW acknowledge funding by Deutsche Forschungsgemeinschaft (DFG, German Research Foundation) – CRC 1456 (project-ID 432680300), project B06. 
AW additionally acknowledges support by DFG RTG 2756 (project-ID 449750155), project A4, and by the German Federal Ministry of Research, Technology, and Space (BMFTR, project PhabiMed, 05M26MGA).}
\thanks{G. Sarnighausen and A. Wald are with the Institute for Numerical and Applied Mathematics, University of Göttingen, Germany.
}
\thanks{A. Hauptmann is with the Research Unit of Mathematical Sciences,
University of Oulu, Oulu, Finland and with the Department of Computer
Science, University College London, London, United Kingdom.
}
\thanks{Codes published at: \url{https://gitlab.com/GesaSarnighausen/transformer-causality-regularization.git}}
\thanks{Manuscript received XXXX, 202X; revised XXXX, 202X.}}



\maketitle

\begin{abstract}
We study the concept of including the causality principle as regularizer into the solution of linear time-dependent inverse problems. This is achieved by combining transformer-based predictions with classical variational regularization, resulting in what we call transformer causality regularization (TCR).
The causality principle states that an object at time $t'$ depends only on its previous states at $t < t'$ and is independent of future states at $t > t'$. Since the transformer architecture represents sequence-to-sequence functions and can be equipped with a causal attention mask, transformers are the natural choice for a learned causality function that predicts the state of an object at time $t'$ given the previous states at $t < t'$. We combine this with the inductive bias of convolutional neural networks (CNNs) for imaging tasks to treat the spatial variable. The output of the spatial-temporal transformer is then used as a prior for variational regularization, such that classical results on regularization and convergence for solution methods directly transfer to our case.
Using the example of dynamic computerized tomography, we compare TCR to a static and dynamic version of the earlier introduced unrolled adversarial regularizer for simulated and measured data. The results show that using TCR within a variational framework improves reconstruction results and data-consistency.
\end{abstract}


\section{Introduction}

Dynamic or time-dependent inverse problems (DIPs) arise in many applications where the quantity of interest evolves over time and is reconstructed from indirect observations. Prominent examples include imaging of process flows \cite{hampel2022review} as well as the estimation of time-dependent parameters in partial differential equations \cite{kaltenbacher17}.

Here, we consider linear dynamic inverse problem, which can be expressed by the operator equation
\begin{align*}
    A(\vartheta) = \psi, \quad A: \X \to \Y,
\end{align*}
where the linear forward operator $A$, the target $\vartheta$ and the data $\psi$ are time- and space-dependent. The spaces $\X$ and $\Y$ are function spaces reflecting the time-dependence of the inverse problem, see, e.g., \cite{Burger_2024, Sarnighausen_2026}. 
In this article, we consider the case of reconstructing the moving object instead of a fixed state. This means that the problem is heavily under-sampled since in general there are only a few measurements available per time step.

Research on dynamic inverse problems is often linked to a concrete application, see~\cite{book_TDPIIP,A1-KLEIN;ET;AL:21,Hauptmann2021}, such as dynamic computerized tomography~\cite{Burger_2017, bh2017, bh2014, luetjen24,hakkarainen2019undersampled}, magnetic resonance
imaging~\cite{hauptmann2019real,kustner2020cinenet}, emission tomography~\cite{hashimoto2019dynamic}, magnetic particle imaging~\cite{albers2023timedependent, A1-KALTENBACHER;ET;AL:21}, or structural health monitoring~\cite{A1-KLEIN;ET;AL:21, C2-LAMBWAVES-BUCH:18,tallman2020structural}; In many cases, additional information on the motion is included to reduce the degrees of freedom in a meaningful way~\cite{bh2014,Burger_2017, Gris_2020, Hahn_2014, bhw20, Nitzsche22, Goedeke_2023}.
We aim to develop methods that do not depend on the specific forward operator and do not assume prior knowledge on the motion.

In recent years, data-driven approaches have become increasingly popular for solving (static) inverse problems~\cite{Arridge_2019}. A variety of approaches have been proposed to incorporate neural networks into the reconstruction process, ranging from purely data-driven methods to hybrid schemes that combine learned components with classical model-based regularization techniques~\cite{Haltmeier2021,mukherjee2023learned}. One possibility is to extend the classical method of variational regularization 
by learning a regularization or penalty term~\cite{Lunz2018, mukherjee2021end, mukherjee2024data, kobler2021total, li2020nett, hauptmann2025convergent}. In classical variational regularization this term needs to be handcrafted and depends on a-priori assumptions on the solution.

In this study, we propose to learn a regularizer using transformer-based architectures. 
Transformers were originally introduced for language processing tasks \cite{17vaswaniAttention} and are the backbone of many foundation models. In image processing tasks, transformers are primarily used for inpainting, denoising, or segmentation~\cite{li2022mat, zhao22_transformerCNNDenoising, 23ZhouNNformer, 21Petit}.
With the growing popularity of transformer-based architectures, transformers have also been used in approaches for solving inverse problems, e.g., for specific static inverse problems such as inverse scattering (post-processing)~\cite{du2024inhomogeneous} or PDE-based inverse problems~\cite{guo2022transformer}.
In~\cite{ahuja2023transformers} the authors solve different static linear inverse problems with in-context learning using transformers. 
Transformers have also been employed for solving dynamic inverse problems. In~\cite{muthukrishnan2023invrt} dynamic radar inverse problems are solved using a spatial-temporal transformer that can capture the temporal dependencies. 
The authors in \cite{sun2025video} trained a spatial-temporal transformer for dynamic scattering media video reconstructions.

The aim of this article is to extend the research on classical solution methods for dynamic inverse problems in Lebesgue-Bochner spaces~\cite{Sarnighausen_2026, Burger_2024}, which are the usual choice for time-dependent partial differential equations, to the deep learning framework and solve dynamic inverse problems using data-driven regularization methods with regularization guarantees. 
In particular, this is achieved by encoding the concept of causality in a spatial-temporal transformer prediction model using causal attention masks. The sequence-to-sequence nature of transformers is particularly suited to deal with a varying length of the temporal information. The transformer predictions are then used to solve a classic variational problem to enforce data-consistency over long time periods. The resulting model is termed \emph{Transformer Causality Regularization} (TCR), see Sections~\ref{sec:varReg} and~\ref{sec:transformerCausalityFunction}. 
We compare our proposed model to the unrolled adversarial regularizer (UAR)~\cite{mukherjee2021end}, a data-driven method developed for static inverse problems, Section~\ref{sec:DL_UAR}, which we extend to a dynamic setting.


To provide a baseline on learned regularization for dynamic inverse problems, we first apply the UAR to a dynamic and static setting in Section~\ref{sec:DL_UAR} without any information on the motion analogously to Tikhonov regularization in Bochner spaces in~\cite{Sarnighausen_2026}. In both cases, classical and data-driven, we understand the object of interest as a time-dependent function. In particular, this means that in the dynamic setting of UAR the measured data is available on the whole spatial-temporal domain. 
In the classical dynamic setting, however, the difference between time and space is expressed by different respective regularities \cite{Sarnighausen_2026}, as is usual in dynamic compressed sensing magnetic resonance imaging \cite{tsao2003k}, i.e., different regularization functionals are chosen for time and space. Nevertheless, in the UAR approach, the difference between time and space is not manually determined in advance but learned by the network. Since the solution is computed for all time points simultaneously in the dynamic setting, the data from later time points can also influence the reconstructions from earlier time points.
This assumption can for some applications contradict the nature of time and the principle of causality, which says that the state of an object $\vartheta \in \X$ at time $t'$, i.e., $\vartheta(t')$, depends on its previous states $\vartheta(t)$ for $t < t'$ only and cannot be influenced by future states. 
If we are, for example, interested in reconstructing a moving object in real-time, it would not be feasible to only get the reconstruction after the movement is finished. Such a causality principle can be implemented through a specialized attention mask, as, e.g., in the spatial-temporal transformer for a specific problem of video reconstruction~\cite{sun2025video}.



This motivates our aim of incorporating the principle of causality into a regularization scheme by learning a causality function that is responsible for including temporal information on the evolution of the dynamics and taking the previously reconstructed time points into account. Nevertheless, the predictions generated purely by this learned causality function only work well for consistent data, i.e., a continuous movement with no change in the dynamical behavior of the object. Therefore, the reconstruction error accumulates over time even for easily predictable motions. Thus, to counteract the accumulating error and to get meaningful predictions even for more complicated motions, we solve a variational problem for each time step ensuring data consistency of the obtained reconstructions.

TCR, our transformer-based approach of variational regularization using a learned prior, includes the principle of causality into the solution of general under-sampled dynamic inverse problems.

We conclude the article by performing numerical experi\-ments for both simulated and real measured experimental data using the example of dynamic computerized tomography, see Sections~\ref{sec:Experiment}-~\ref{sec:discussion}.

\section{Dynamic Inverse Problems}
\label{sec:DLMethodsDIP}

\subsection{Causality regularization}
\label{sec:varReg}

We aim to include the principle of causality into a data-driven regularization method for dynamic inverse problems, i.e., the state of an object $\vartheta \in \X$ at time $t'$, i.e. $\vartheta(t')$, depends on its previous states $\vartheta(t)$ for $t < t'$ and cannot be influenced by future states.

Therefore, we want to find a suitable causality function $\Lambda$ with
\begin{align*}
    \Lambda \left( t', \{\vartheta(t): 0 \leq t < t' \} \right) \approx \vartheta(t'),
\end{align*}
which predicts the state at time $t'$ given the previous states at times $t < t'$. This function is then included as a prior in a penalty term in a classical variational regularization method, which consists in solving the minimization problem
\begin{align}
\label{eq:DIP-NN-minProblem}
\begin{split}
    \min_{\vartheta(t')} \Bigl\{ & D\left[ A_{t'} \vartheta(t') - \psi(t') \right] \\
    &+ \alpha_{t'} P\left[ \vartheta(t') - \Lambda\left(\{\vartheta(t): 0 \leq t < t' \} \right) \right] \Bigr\},
\end{split}
\end{align}
with functions $D$ and $P$ for the data discrepancy and penalty term, respectively, a regularization parameter $\alpha_{t'}$ for all $t' \in (t^*,\tilde{T}]$ with $t^* > 0$, $A_{t'}$ the forward operator at time $t'$, $\psi(t')$ the measured data at time $t'$ and $\vartheta(t')$ the function we want to recover at time $t'$. This approach is built on the availability of sufficient information on initial reconstructions in a time interval $[0,t^*)$. For the specific transformer causality function from Section~\ref{sec:transformerCausalityFunction} we call the proposed method transformer causality regularization (TCR).


\subsubsection{Minimizing the variational problem}
\label{sec:minVarProb}
The minimization of~\eqref{eq:DIP-NN-minProblem} depends on the choice of $D$ and $P$, i.e., the data discrepancy and penalty term, and can be achieved analogously to the standard case of variational regularization for $\Lambda = 0$ since $\Lambda$ does not depend on $\vartheta(t')$. Thus regularization and convergence results for the classical case with $\Lambda = 0$ directly transfer to causality regularization. In particular, causality regularization is a regularization method for each time point $t'$. For clarity and completeness, we reformulate the classical algorithms employed here and describe in detail how the causality function is incorporated.
We mainly consider $L^2$- and $L^1$-regularization, i.e., we choose
\begin{align}
\label{eq:L2reg}
    D(\cdot) &= P(\cdot) = \frac{1}{2} \| \cdot \|^2_2, \quad \text{or} \\
    \label{eq:L1reg}
    D(\cdot) &= \frac{1}{2} \| \cdot \|^2_2, \quad P(\cdot) = \| \cdot \|_1.
\end{align}

For the choice~\eqref{eq:L2reg}, we minimize the Tikhonov functional~\eqref{eq:DIP-NN-minProblem} by gradient descent. Since $\Lambda$ does not depend on $\vartheta(t')$, the gradient of~\eqref{eq:DIP-NN-minProblem} is given by
\begin{align*}
    &(A_{t'})^* \left(A_{t'} \vartheta(t') - \psi(t') \right) \\
    &+ \alpha \Bigl(\vartheta(t') - \Lambda\left(\{\vartheta(t): 0 \leq t < t' \} \right) \Bigr).
\end{align*}
For this case, we propose the following algorithm:

\begin{algo}{($L^2$-causality regularization)}
\label{algo_L2-TCR}
Given $t' \in (t^*,\tilde{T}]$, previous reconstructions $\vartheta(t)$ for $t \in [0, t^*]$, regularization parameters $\alpha_t$ and initial guesses $\vartheta(t')_0$, the reconstruction $\vartheta(t')$ can be computed by the iteration
\begin{align}
    \label{eq:TikhonovNNIteration}
    \begin{split}
    \vartheta(t')_{k+1} = \vartheta(t')_k &- \tau_t \biggl( (A_{t'})^* (A_{t'} \vartheta(t')_k - \psi(t')) \\
    &+ \alpha_{t'} \Bigl(\vartheta(t')_{k} - \Lambda\left(\{\vartheta(t): 0 \leq t < t' \} \right) \Bigr) \biggr),
    \end{split}
\end{align}
where the step size $\tau_t$ is chosen as $\tau_t = \frac{1}{\lVert A_t^* A_t \rVert}$ in the operator norm to guarantee convergence as in the standard Landweber iteration scheme, see e.g.~\cite{Rieder:keineProblemeMitInversenProblemen, Engl1996RegularizationOI}.
\end{algo}

For the $L^1$-penalty term~\eqref{eq:L1reg}, we can minimize~\eqref{eq:DIP-NN-minProblem} using an accelerated proximal algorithm, introduced in~\cite{FISTA,daubechies2004iterative}.

\begin{algo}{($L^1$-causality regularization)}
\label{algo:FISTA}
   Given $t' \in (t^*,\tilde{T}]$, previous reconstructions $\vartheta(t)$ for $t \in [0, t^*]$, regularization parameters $\alpha_t$, the continuous, convex and non-smooth functional $G:\R^n \to \R$ with $G(x) = \alpha_t \| x -  \Lambda\left(\{\vartheta(t): 0 \leq t < t' \} \right) \|_1$ and the differentiable functional $F: \R^n \to \R$ with $F(x) = \frac{1}{2} \| A_{t'} x - \psi(t')\|^2_2$ with Lipschitz continous gradient $\nabla F$, we want to minimize
    \begin{align}
        \label{eq:fistaMIN}
        \min_{x \in \R^n} \Bigl\{ \frac{1}{2} \| A_{t'} x &- \psi(t')\|^2_2\Bigr. + 
        \\        
        \Bigl.&\alpha_t \| x - \Lambda\left(\{\vartheta(t): 0 \leq t < t' \} \right) \|_1 \Bigr\},
        \nonumber 
    \end{align}
    provided a solution to~\eqref{eq:fistaMIN} exists. 
    This can be done using the proximal operator defined by $\text{prox}_{G}\bigl(v) = \argmin_x \left( G(x) + \frac{1}{2} \| x - v \|^2_2 \right)$ with the following iteration scheme:

\begin{figure}[H]
\centering
\begin{algorithmic}[1]
\STATE \textbf{Input:} $l$ Lipschitz constant of $\nabla F$ \\
\STATE $y_1 \gets x_0$, $t_1\gets 1$ \\
\STATE choose step size $\lambda \in (0, 1/l]$ for convergence \\
\FOR{$k \geq 1$}
    \STATE $x_k \gets \text{prox}_{\lambda G}\bigl(y_k - \lambda \nabla f (y_k)\bigr)$ \\
    \STATE $h_{k+1}  \gets \frac{1 + \sqrt{1 + 4 h_{k}^2}}{2}$ \\
    \STATE $y_{k+1} = x_k + \frac{t_{k} - 1}{h_{k+1}} (x_k - x_{k-1})$
\ENDFOR
\STATE \textbf{return:} $x_k$
\end{algorithmic}
\end{figure}
\end{algo}

We choose $l = \| A_{t'}^* A_{t'} \|$ with the operator norm for convergence. The proximal operator $\text{prox}_{\lambda \| \cdot - \Lambda\left(\{\vartheta(t): 0 \leq t < t' \} \right)  \|_1}$ can be easily computed using the fact that the proximal operator of a shifted function $H = G(\cdot - a)$ is given by $\text{prox}_H(v) = a + \text{prox}_G(v -a)$, see, e.g.,~\cite[Proposition 24.8]{Bauschke2017Proximity}, and the proximal operator of the $L^1$-norm is the soft thresholding operator \cite[Chapter 6.5]{14ParikhProximalAlgos}.

For the application to experimentally measured data in Section~\ref{sec:Experiments_rollingStonesData}, we add an additional $TV$-term to the minimization problem, i.e., we solve
\begin{align}
\label{eq:TVL1Min}
\begin{split}
     \min_{\vartheta(t')} \Bigl\{ &\left \| A_{t'} \vartheta(t') - \psi(t') \right \|_2^2 + \beta_{t'} \| \nabla \vartheta(t') \|_{1} \\
     &+ \alpha_{t'} \left \| \vartheta(t') - \Lambda\left(\{\vartheta(t): 0 \leq t < t' \right) \right \|_1    \Bigr\}.
\end{split}     
\end{align}

Similarly to $L^1$-causality regularization we use a proximal algorithm to solve~\eqref{eq:TVL1Min}; in this case the primal-dual hybrid gradient method~\cite{2011ChambollPDHGAlgo} for minimizing
 \begin{align}
        \label{eq:PDHGMIN}
        \min_{x \in \R^n} F(Lx) + G(x),
    \end{align}
    where $G:X \to [0,\infty]$, $F: Y \to [0,\infty]$ are convex and lower-semi-continuous functionals, $L : X \to Y$ is a continuous linear operator, and the spaces $X$ and $Y$ are finite dimensional. 

To obtain a $L^1$-TV-causality regularized solution to~\eqref{eq:TVL1Min}, we choose $L = (I, \nabla): X \to X \times X^d$ for a $d$-dimensional space $X$ and the identity $I$ on $X$, $F(x,u) = \frac{1}{2} \| x - \psi(t') \|^2_2 + \beta_t \| u \|_1$ and $G(x) = \alpha_t \| x -  \Lambda\left(\{\vartheta(t): 0 \leq t < t' \} \right) \|_1$.

This choice for $L$ and $F$ is standard for $TV$-regularization. Details can be found in~\cite{2011ChambollPDHGAlgo} and in the documentation for the operator discretization library 'odl'~\url{https://odlgroup.github.io/odl/math/solvers/nonsmooth/pdhg.html} that we used for solving the problem. The additional $L^1$-term in the functional $G$ can then be treated as in the $L^1$-case. 

\subsection{Unrolled adversarial regularizer (UAR) for the dynamic case}
\label{sec:DL_UAR}
The idea of UAR is to combine iterative unrolling with data-adaptive regularization via a generative adversarial framework as proposed in~\cite{2021Mukherjee}. The UAR approach extends the original adversarial regularizer~\cite{Lunz2018} which only consisted of the regularizer (discriminator).

The concept of generative adversarial networks (GANs) is to train two networks alternately in an unsupervised way: One generator (here: network for computing reconstructions) and one discriminator (here: regularizer). While the generator aims to generate 'good' samples from a distribution, the discriminator is trained to discriminate between real samples of the distribution and generated ones by the generator. The discriminator is optimized using the output of the generator, whereas the generator is optimized taking into account the assessment of the discriminator. Since this scheme is iterated, the discriminator gets better at distinguishing real and generated samples and the generator produces samples closer to the real distribution. In theory, this continues until the output of the generator is indistinguishable from the real distribution.

The generator $G_\genPar$ is implemented as an iterative unrolling strategy that aims to minimize the variational problem via a primal-dual type method using two CNNs $\Gamma^{\text{primal}}_{\genPar}$ and $\Gamma^{\text{dual}}_{\genPar}$ for the primal and dual optimization each, where $\genPar$ are the network parameters. The networks $\Gamma^{\text{primal}}_{\genPar}$ and $\Gamma^{\text{dual}}_{\genPar}$ correspond to learned versions of proximal operators, see Section~\ref{sec:minVarProb} for a classical primal dual algorithm. For $L = 20$ layers (iterations) the next iterates for $1 \leq l \leq L-1$ are computed as
\begin{align*}
    h_{l+1} &= \Gamma^{\text{dual}}_{\genPar} \left( h_l, \sigma_l, A(\vartheta_l), \psi^\delta \right), \\
    \vartheta_{l+1} &= \Gamma^{\text{primal}}_{\genPar}  \Bigl( \vartheta_l, \tau_l, A^*(h_l) \Bigr), 
\end{align*}
where we choose $h_0 = 0$ and $\vartheta_0 = A^+ \psi^\delta$ for the pseudo inverse $A^+$.
The step sizes $\sigma_l$ and $\tau_l$ are trainable parameters both initialized with $0.01$. 

The discriminator or regularization network $R_\regPar$ is implemented as a CNN with parameters $\regPar$ and $6$ convolutional layers followed by average pooling and two dense layers.

In the $(k+1)$-st iteration of UAR the parameters of the regularizer $R_\regPar$ are optimized for fixed $\genPar_{k}$ as
\begin{align} 
\label{eq:minRegularizer}
    \regPar_{k+1} &\in \argmin_{\regPar : R_{\regPar} \in \mathbb{L}^1}  \{ J_{k}^{\text{reg}} (\regPar )\} , \quad \text{with}\\
     J^{\text{reg}}_{k} (\regPar ) &= \mathbb{E}_{\pi_{\psi^\delta}} \left[ R_\regPar (G_{\genPar_{k}}(\psi^\delta)) \right] - \mathbb{E}_{\pi_{\vartheta}} \left[ R_\regPar (\vartheta) \right],
     \label{eq:lossRegularizer}
\end{align}
where $\mathbb{L}^1$ denotes the space of $1$-Lipschitz functions, $\pi_{\psi^\delta}$ the distribution of noisy measurements, $\pi_\vartheta$ the ground truth distribution and $\mathbb{E}_{\pi_{\psi^\delta}} $ is the expected value. The parameters of the regularization operator $G_{\genPar}$ are afterwards optimized for fixed $\regPar_{k+1}$ as
\begin{align} 
\label{eq:minGenerator}
    \genPar_{k+1} &\in \argmin_{\genPar}  \{ J_{k+1}^{\text{gen}} (\genPar )\} , \quad \text{with}\\
     J^{\text{gen}}_{k+1} (\genPar ) &= \mathbb{E}_{\pi_{\psi^\delta}} \left[ \| \psi^\delta - A(G_{\genPar}(\psi^\delta)) \|^2_2 + \alpha R_{\regPar_{k+1}} (G_\genPar(\psi^\delta)) \right].
     \label{eq:lossGeneraor}
\end{align}
 To compute the loss functionals~\eqref{eq:lossRegularizer} and~\eqref{eq:lossGeneraor} it is not necessary that data pairs are matched. Thus the training works in an unsupervised manner. 

In the actual training, we use discretized versions of~\eqref{eq:lossRegularizer} and~\eqref{eq:lossGeneraor} for computing the loss
            \begin{align}
            \label{eq:lossRegDiscr}
                \tilde{J}^{\text{reg}}(\regPar) &= R_{\regPar}(\vartheta) - R_{\regPar}(u) + \lambda_{gp} \left( \| \nabla R_{\regPar}(\vartheta^\varepsilon) \|_2 - 1 \right)^2 \\
                            \label{eq:lossGenDiscr}
                \tilde{J}^{\text{gen}}(\regPar) &= 
                \bigl \| \psi - A\left( G_{\genPar}(\psi) \right) \bigr\|^2_2
                + \alpha R_\regPar \left( G_{\genPar} (\psi) \right).
            \end{align}
For computational reasons, we only refine the network parameters $\genPar$ and $\regPar$ in each adversarial iteration by performing one optimizer update, instead of fully optimizing each functional for every time step. Further details on the training procedure can be found in~\cite{mukherjee2021end}.

\section{Transformer causality function}
\label{sec:transformerCausalityFunction}
In the following, we consider a time-discretized setting with $T$ equidistantly spaced time steps in a time interval $[0,\tilde{T}]$. The forward operator, measured data and function we want to recover at time $\frac{\tilde{T} t}{T-1}$ for $t = 0 , \ldots , T-1$ are denoted by $A_t$, $\psi^t$ and $\vartheta^t$ respectively. Furthermore, we assume that we already have initial reconstructions for the first two time points $\hat{\vartheta}^0$ and $\hat{\vartheta}^1$. Then our model $\Lambda_{\text{pre}}$ predicts the next time point given the previous reconstructions.

Since our algorithm requires two initial reconstructions of possibly poor quality, we additionally want to learn a model $\Lambda_{\text{re}}$ that produces refined versions of the initial reconstructions which then serve as input to the prediction model during training.

Altogether, we arrive at Algorithm~\ref{alg:VarRegLearnedPrior}:

\begin{algo}{(Transformer causality regularization - TCR)}
\begin{figure}[H]
\centering
\begin{algorithmic}[1]
\STATE \textbf{Input:} $\hat{\vartheta}^0, \hat{\vartheta}^1$ initial reconstructions for the first two time steps, $\Lambda_{\text{re}}, \Lambda_{\text{pre}}$ trained models for refinement of the initial reconstructions and predicting the next frame \\
\STATE \textit{a) compute refinement of initial reconstructions} \\
\STATE $\vartheta^0_{\text{re}}, \vartheta^1_{\text{re}} \gets \Lambda_{\text{re}}(\hat{\vartheta}^0, \hat{\vartheta}^1)$\\
\FOR{$t \in \{0,1 \}$}
    \STATE  $ \vartheta^t \gets \min_{f} \left\{ D(A_{t} f - \psi^t) + \alpha_t P( f - \vartheta^t_{\text{re}} ) \right\}, $
\ENDFOR
\STATE \textit{b) compute reconstructions for previous time steps}\\
\FOR{$t \in \{ 1, \ldots, T-1 \}$}
       \STATE $ \vartheta^t_{\text{pre}} \gets \Lambda_{\text{pre}}(\vartheta^0, \ldots, \vartheta^{t-1})$  
    
      \STATE $ \vartheta^t \gets \min_{f} \left\{ D( A_{t} f - \psi^t ) + \alpha_{t} P( f - \vartheta^t_{\text{pre}} )  \right\}, $
\ENDFOR

\STATE \textbf{return:} $\vartheta^t$ for $t \in \{0,\ldots,T-1\}$
\end{algorithmic}
\end{figure}
\label{alg:VarRegLearnedPrior}
\end{algo}
For data discrepancy $D$ and penalty $P$ as in~\eqref{eq:L2reg} and~\eqref{eq:L1reg} we then solve the variational problem with 
Alg. \ref{algo_L2-TCR} and Alg. \ref{algo:FISTA} to obtain the final versions of the algorithms as $L^2$-TCR and $L^1$-TCR, respectively.

\paragraph{Architecture of spatial-temporal transformer}
\label{sec:TransformerArchitecture}

\begin{figure*}[htb]
\centering

\subfloat{
\includegraphics[width=0.34\textwidth]{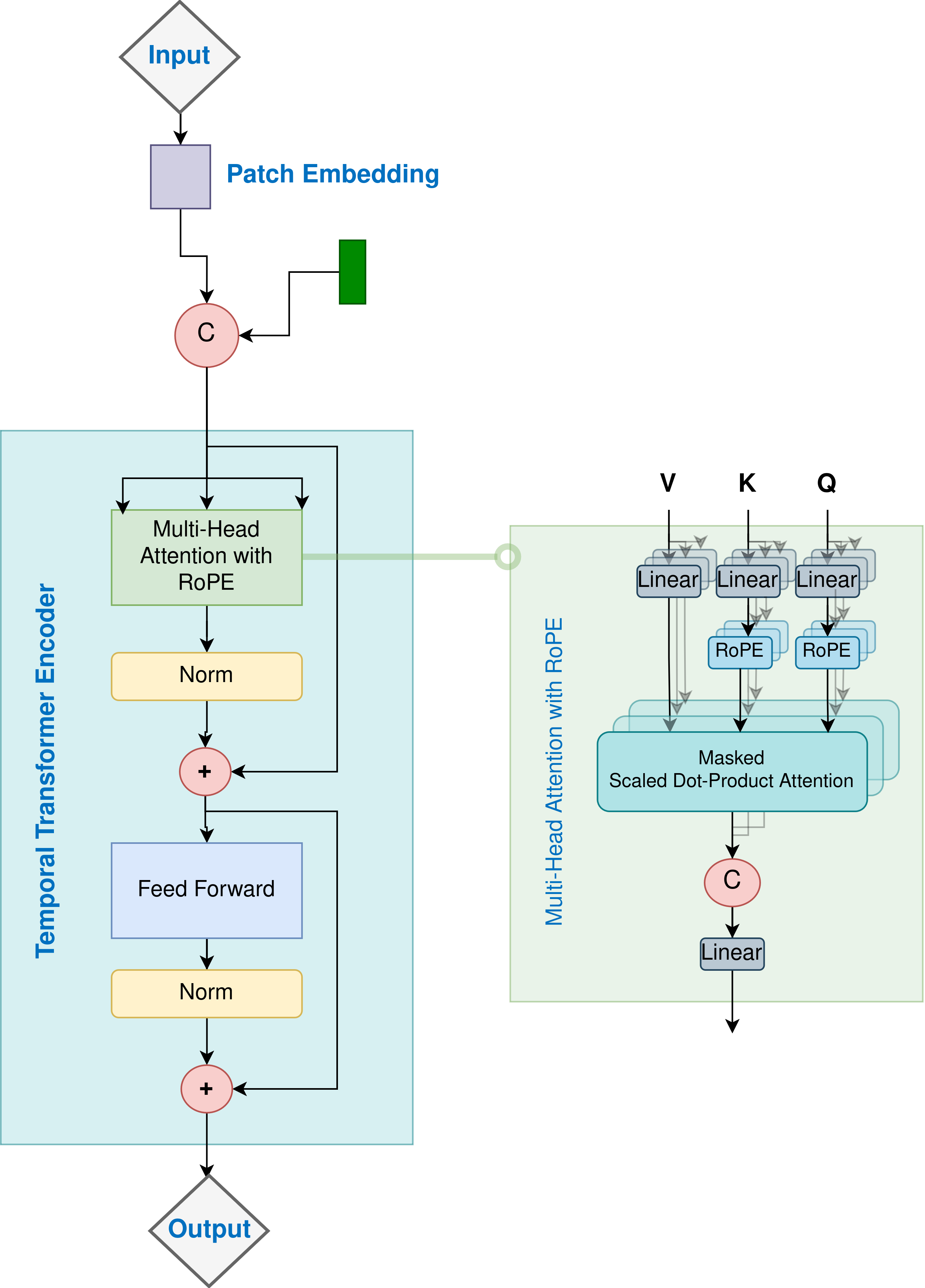}
\label{fig:transformer}
}
\hfill
\subfloat{
\includegraphics[width=0.63\textwidth]{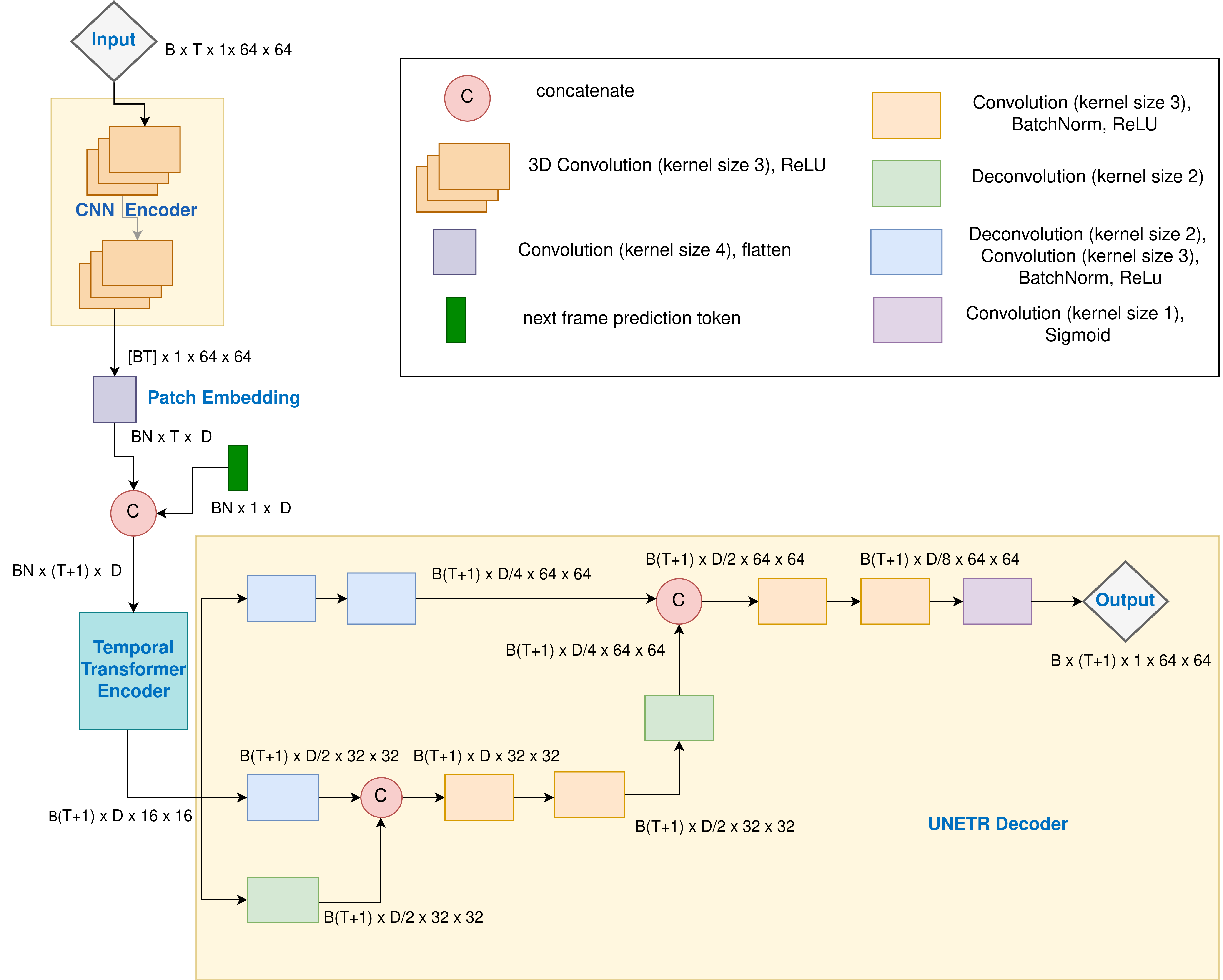}
\label{fig:architecture}
}

\caption{ \textbf{Left:} Architecture of temporal transformer encoder with rotary position embedding (RoPE) with $h = 8$ heads and $6$ transformer layers. \\
\textbf{Right:} Architecture of next frame prediction/refinement model, Number of patches N=356, transformer model dimension D = 512, T length of input sequence, B = 8 batch size for training}
\label{fig:model_architecture}

\end{figure*}

Our architecture aims to combine the inductive bias of CNNs for imaging tasks with the causality of transformers. This allows treating the spatial and temporal dimension differently and naturally transfers the Bochner space setting discussed in~\cite{Sarnighausen_2026} to data-driven methods.
In particular, we use a CNN encoder and a CNN-based UNETR decoder similar to~\cite{UNETR} in the beginning and at the end to deal with the spatial properties of our input sequence, see Figure~\ref{fig:architecture}. The CNN encoder consists of 3D convolutions and thus already preprocesses the temporal dimension additionally to the spatial dimension.  

After a patch embedding that maps the outputs of the CNN encoder to the transformer model dimension $D = 512$, the spatial components are moved to the batch dimension such that the temporal transformer~\ref{fig:transformer} can attend to the temporal properties of the sequence. The transformer encoder is based on the original architecture from~\cite{17vaswaniAttention}. Instead of using additive positional encodings, we make use of rotary position embeddings (RoPE) as in~\cite{24RoPE}, which are applied in the multi-head attention block, see Figure~\ref{fig:transformer}, and combine positional information on the absolute and relative position of each token. As in~\cite{17vaswaniAttention} we use 6 identical transformer layers after each other that are identical to the one depicted in~\ref{fig:transformer}.

The idea of combining CNNs and transformers is commonly employed.
For instance, in~\cite{23ZhouNNformer} a transformer for volumetric medical image segmentation is presented that exploits the combination of interleaved convolution and self-attention operations, and in \cite{21Petit} the authors make use of an overall U-shaped architecture with attention mechanisms. CNNs and transformers are also combined for image denoising. \cite{TIAN2024_CrossTransformerDenoising} proposes a network architecture that employs CNNs and transformer mechanisms alternately. In \cite{zhao22_transformerCNNDenoising} a Swin (Shifted Window) Transformer encoder is combined with a CNN decoder. However, instead of interleaving convolution and self-attention as in the discussed literature, we mainly use convolutions for the spatial and attention for the temporal variable.

The advantage of a transformer for temporal tasks is that the input length is not fixed in contrast to most neural network architectures since all inputs are mapped to a fixed model dimension $D$ of the transformer. This allows to use all previous reconstructions as input for the transformer and not only a fixed amount.  In~\cite{20YunTransformerApproximator} it was shown that transformers are universal approximators of continuous permutation equivariant sequence-to-sequence functions with respect to the Lebesgue-Bochner metric
\begin{align}
\label{eq:metricApprox}
    d(\bar{\vartheta},\bar{\psi}) = \left(\int \| \bar{\vartheta(t)} - \bar{\psi(t)} \|_p^p \dt \right)^{\frac{1}{p}}
\end{align}
for $\bar{\vartheta}, \bar{\psi} \in \R^{n_t \times n_d}$ with $n_t$ the temporal and $n_d$ the spatial size, which transfers directly to $\R^{n_t \times n_d \times n_d}$. Note that the integral in~\eqref{eq:metricApprox} for discrete $\bar{\vartheta}, \bar{\psi}$ can be understood as a sum. As discussed in~\cite{Burger_2024, Sarnighausen_2026}, the Lebesgue-Bochner setting is natural for time-dependent problems, making transformers an obvious choice for our application. 

However, the Lebesgue-Bochner setting does not include positional information or the principle of causality since the object is seen as an $(N+1)$-dimensional object, which is reconstructed for all time steps simultaneously, for $N$ spatial dimensions. In the transformer case, this is reflected by the class of permutation equivariant functions.

To include positional information in the transformer architecture, positional encodings are typically added~\cite{17vaswaniAttention, 24RoPE}.
Together with a causal mask used during training, this only allows the model to attend to previous time-steps and aligns with the causality principle. In~\cite{sun2025video} causality is also enforced by a specialized attention mask.
The authors of~\cite{20YunTransformerApproximator} showed that transformers equipped with a learnable additive positional encoding are universal approximators of continuous functions that map a compact domain of $\R^{n_t \times n_d}$ to $\R^{n_t \times n_d}$. Even though we are not using a learnable additive positional encoding, but a multiplicative one, this property motivates the use of transformers for our application.
The universality results of~\cite{20YunTransformerApproximator} concern shallow transformers with $2$ heads and thus require that the transformer's embedding dimension grows with the number of tokens. The authors of~\cite{furuya2025transformers} take a different focus and show that deep transformers with a fixed embedding dimension are universal for an arbitrary large number of tokens.

Since we separate the use of CNNs and transformers, in principal the approximation results can be transferred to the transformer part of our setting. In particular, we use the CNN encoder to obtain the temporal tokens that serve as input to the transformer. The transformer's output is then post-processed by the CNN based decoder. Thus, we apply the transformer only as a temporal sequence-to-sequence function.

\paragraph{Training schemes}

In the following, we train neural networks for the functions $\Lambda_{\text{re}}$ and $\Lambda_{\text{pre}}$. 
 For simplicity, we use the same network architecture, the spatial-temporal transformer described in the previous Section~\ref{sec:TransformerArchitecture}, see Figure~\ref{fig:architecture}, for both the refinement and the prediction model.
 In both cases the input consists of reconstructions up to time $t$ and the output are the reconstructions up to time $t+1$, i.e.,
 \begin{align*}
    \Lambda_{\{\text{re}, \text{pre} \}} : \R^{T \times  H \times W} \to \R^{(T+1) \times H \times W}.
 \end{align*}
 Only the loss and training routine change for both models.

 For the refinement model, we are purely interested in the refinement of the first two time steps. Therefore, we choose the loss as
\begin{align}
\label{eq:loss_refine}
\text{loss}_{\text{re}} = \frac{1}{2} \Bigl( \| \vartheta^0_{\text{re}} - \vartheta^0_{\text{gt}} \|_{2}^2 + \| \vartheta^1_{\text{re}} - \vartheta^1_{\text{gt}} \|_{2}^2
\Bigr),
\end{align}
where $\vartheta^t_{\text{gt}}$ is the ground truth at time point $t$ for $t=0,1$, and the input sequence always consists of the first two time-steps, either as ground truth or as noisy reconstructions obtained by classical methods, such that the loss is computed only with respect to the ground truth data of time steps $0$ and $1$.

In the prediction model we are not interested in the output of the previous time steps but only in the prediction for the next time step, i.e., we choose the loss as
\begin{align*}
    \text{loss}_{\text{pre}} = \frac{1}{r\_s} \sum_{t = 2}^{r\_s + 1} \| \vartheta^t_{\text{pre}} - \vartheta^t_{\text{gt}} \|_{2}^2 ,
\end{align*}
where the rollout steps $r\_s$ determine the number of prediction steps. For example, for two rollout steps we get the first prediction for $t=2$ and use this (or the ground truth) to predict the reconstruction at $t=3$. To reduce computational complexity, we optimize the model parameters after each rollout step.

\section{Experiments}
\label{sec:Experiment}

In the following, we apply the methods described in Section~\ref{sec:DLMethodsDIP} to the inverse problem of dynamic computerized tomography.

\subsection{Model of dynamic computerized tomography}
The static Radon transform $R$, the forward operator of computerized tomography, integrates a function over all affine hyperplanes (lines in 2D), and is defined by
\begin{align*}
(R f)(\varphi, \sigma) = g(\varphi, \sigma) = \int_{-1}^1  f\left(\sigma \theta(\varphi) + w \theta(\varphi)^{\bot}\right) \dw
\end{align*}
 for $1 < r < \infty$, $\varphi \in [0, \pi)$, $\sigma \in [-1,1]$, the unit circle $B \subset \R^2$ with $\theta(\varphi) = ( \cos(\varphi), \sin(\varphi))$ and $\theta(\varphi)^{\bot} = ( -\sin(\varphi), \cos(\varphi))$. This choice is not a constraint since each bounded domain in $\R^2$ can be scaled to the unit circle. 
Now, we consider the dynamic Radon transform defined by
\begin{align*}
    (\dynR  \vartheta)(t)  = \psi(t) \coloneqq R (\vartheta(t)) 
\end{align*}
for each $t \in [0,\tilde{T}]$.

In the discretized case, we only have a few angles per time step available. In the following, we consider $100$ equidistantly spaced offset values in $[-1,1]$ discretizing the sensor and either $3$ or $10$ equidistantly spaced and rotating angles per time step, such that the projection angles are different for each time step.
To use TCR (Algorithm~\ref{alg:VarRegLearnedPrior}), we additionally require initial reconstructions of the first two time steps with sufficient quality. To this end, we assume that 20 projections are available for each of these two time steps. While this choice permits sufficient reconstruction quality for our purpose, we would like to emphasize that this number of projections still leads to sparse data.



\subsection{Training data}

We use the same dynamic phantom data consisting of 10 time steps and $64 \times 64$ pixels for both TCR and UAR. The underlying motion is a random affine linear motion, and each phantom contains 3-5 rectangles, ellipses, and/or circles. The training data set consists of $5000$ phantoms and the test data set of $50$ phantoms.

In the unsupervised UAR approach, we use unpaired sinograms and ground truth as training data sets to sample from the distribution of noisy measurements (sinograms) and ground truth, whereas for supervised TCR, we use noisy Landweber reconstructions obtained from $20$ measured angles and ground truth data for the first two time steps (depending on teacher forcing, see~\eqref{eq:teacher_forcing}) as input and ground truth as output; see Figure~\ref{fig:Phantom4} as an example.

\subsection{Training details}

\subsubsection{Models for UAR}
\label{sec:TrainingDetailsUAR}
Similarly as in the purely analytical approach in~\cite{Sarnighausen_2026}, we consider a dynamic setting where a solution to the dynamic inverse problem is computed for all time points simultaneously, and a static setting where each time point is considered independently.

For all three CNNs from the networks $G_\genPar$ and $R_\regPar$ we use $3D$ convolutions in the dynamic setting where the input contains sequences of images and $2D$ convolutions in the static setting where the input consists only of individual images. 

The dynamic data set consists of the entire time sequence for ground truth and measurement data, whereas for the static data set, we only take one randomly selected time step from each pair of ground truth and measurement data.

In the training procedure, we choose a batch size of 1, meaning that we optimize for one sample only so that we can ignore the expected value in~\eqref{eq:lossRegularizer} and~\eqref{eq:lossGeneraor}. We first train an initial baseline regularizer and reconstruction operator by optimizing the respective loss functions~\eqref{eq:lossRegDiscr} and~\eqref{eq:lossGenDiscr} for 5 epochs each using Adam with learning rate $\eta = 10^{-5}$ and optimizer parameters $(\beta_1,\beta_2) = (0.5,0.99)$. In the following adversarial training, we train for 10 epochs and update the parameters of $R_\regPar$ and $G_\genPar$ alternately with learning rate $\eta = 2\cdot 10^{-5}$ and $(\beta_1,\beta_2)$ as before.
For further details on the training procedure we refer to~\cite{mukherjee2021end}.

\subsubsection{Models for TCR}

Since we want to use the architecture of the spatial temporal transformer described in Section~\ref{sec:TransformerArchitecture} for two different tasks, namely refining the first two reconstructions and predicting the next reconstruction for a given input sequence of varying length, we use a different loss function and training procedure for both tasks:

\paragraph{Training of the refinement model}

We train the model for 100 epochs in total, starting with ground truth images and then gradually moving to Landweber reconstructions computed from $20$ equidistantly spaced and rotating angles and $100$ offset values. The ground truth ratio (gt\_r) determines the proportion of ground truth data in the training data for one epoch and is defined by
\begin{align*}
\text{gt\_r}(\text{e}) = 
    \begin{cases} 
      1 & \text{e} <  10 \\
      1.0 -  \frac{\text{e} - 10}{37.5} & 10\leq \text{e} < 40 \\
      0.2 & \text{e} \geq 40,
   \end{cases}
\end{align*}
for each epoch e.

For the learning rate, we choose a cosine annealing schedule as proposed in~\cite{16SGDR_cosineAnnealing} without restarts, 10 warm up epochs, 100 epochs in total, $\text{min\_lr} = 1e\!-\!6$ and $\text{max\_lr} = 1e\!-\!4$.
As an optimizer, we use the AdamW algorithm as proposed in~\cite{19adamw} with $\beta=(0.9, 0.95)$ and weight decay$=1e\!-\!2$.

\paragraph{Training of the prediction model}

To train the prediction model, we use the output of the refinement model, i.e., the refined reconstructions for the first two time steps, as input for the prediction model. To stably train the model auto-regressively, we use probabilistic teacher forcing, i.e., we start predicting the next frame for time steps $t=\{2,\ldots,r\_s +1\}$ from ground truth inputs, gradually decrease the proportion of ground truth inputs and use the predicted outputs instead. Additionally, we increase the number of maximum rollout steps and the probability that rollout is actually happening. If not, we only compute the prediction for $t=2$. This accelerates training and ensures that the task is not too hard in the beginning.

The maximum number of rollout steps is set to 2 until epoch 30, increased to 4 until epoch 70, to 6 until epoch 90, and to 8 for the remaining epochs.
The maximum number of possible rollout steps is 8 for our dataset. The rollout probability increases linearly to 1 over 30 epochs and can be computed by $\min(1, \text{e} / 30) $. The teacher forcing ratio (tf\_r) decreases linearly from $0.9$ to $0$ until epoch 85 via 
\begin{align}
\label{eq:teacher_forcing}
    \text{tf\_r}(\text{e}) = 
    \begin{cases}
        \max\{0, 0.9 (1 - \text{e}/85)\} & \text{e} < 85 \\
        0 & \text{else}
    \end{cases}.
\end{align}
Regarding the refinement model, we train for $100$ epochs in total and use the AdamW algorithm~\cite{19adamw} as an optimizer with $\beta=(0.9, 0.95)$ and weight\_decay$=1e\!-\!2$.
Here, we use a cosine annealing scheduler that starts at $3e\!-\!5$ and only starts decreasing after $40$ epochs to $1e\!-\!6$.

\paragraph{Choosing the regularization parameters}

To apply TCR to the test data set, we need to choose regularization parameters $\alpha_t$ for every time step. To ensure consistency over time, we only choose two different parameters: one for the first time steps and one for the rest, since we have more data available for the first two time steps. We select the parameters that minimize the $L^2$-error with respect to the ground truth. 
For $L^1$-TCR, we stop the algorithm after 200 iterations, for $L^2$-TCR after a maximum of 19 iterations if $\| A_t \vartheta^{t}_k - \psi^t \| > \| A_t \vartheta^{t}_{k-1} - \psi^t \| $ is not fulfilled before, i.e., the data discrepancy starts to increase.

\subsection{Implementation}
For the implementation of both UAR and TCR, we use the operator discretization library 'ODL'~\cite{jonas_adler_ODL} which in particular contains implementations of FISTA and PDHG.

For the implementation of UAR we use the existing codes of the authors of~\cite{2021Mukherjee}, publicly available at~\url{https://github.com/Subhadip-1/unrolling_meets_data_driven_regularization}, as a basis which we extend to a dynamic setting and our data set.

\subsection{Reconstruction from experimental data}
\label{sec:Experiments_rollingStonesData}
We apply TCR to the 'rolling stones data set' presented in~\cite{Burger_2017}, which consists of three stones drifting away from each other over 30 time steps, see Figure~\ref{fig:rollingStones}. The measurements are taken in a $\mu$CBCT system. The flat stones are centered right in front and the middle of the flat panel detector. This results in a source-to-detector distance much larger than object-to-detector distance and hence can be well approximated by a parallel beam geometry in 2D. 
We use sinograms with 63 offset values spaced in [-1,1], 10 equidistantly spaced angles for the first two time steps and 3 angles for the remaining time steps, respectively. To choose the regularization parameters for $t= 0,1$ and $t > 1$, we first manually select parameters for the initial reconstruction and then use this to hand tune parameters for the other time steps. The reconstructions are of size $64 \times 64$ to ensure compatibility with the transformer causality function.

We note, that the motion of this dataset does not follow the training data, as all three stones are individually moving. Additionally, the measurements are taken over 30 times steps instead of only 10. Thus, this is an ideal dataset to test generalizability of our proposed method to experimentally measured out-of-distribution data.

\section{Results}
\label{sec:result}

\subsection{Reconstruction from simulated data}

\subsubsection{UAR}
\begin{figure*}[htbp]
    \centering
 {
    
        \includegraphics[width=\textwidth,  trim=300 230 300 200,
  clip]{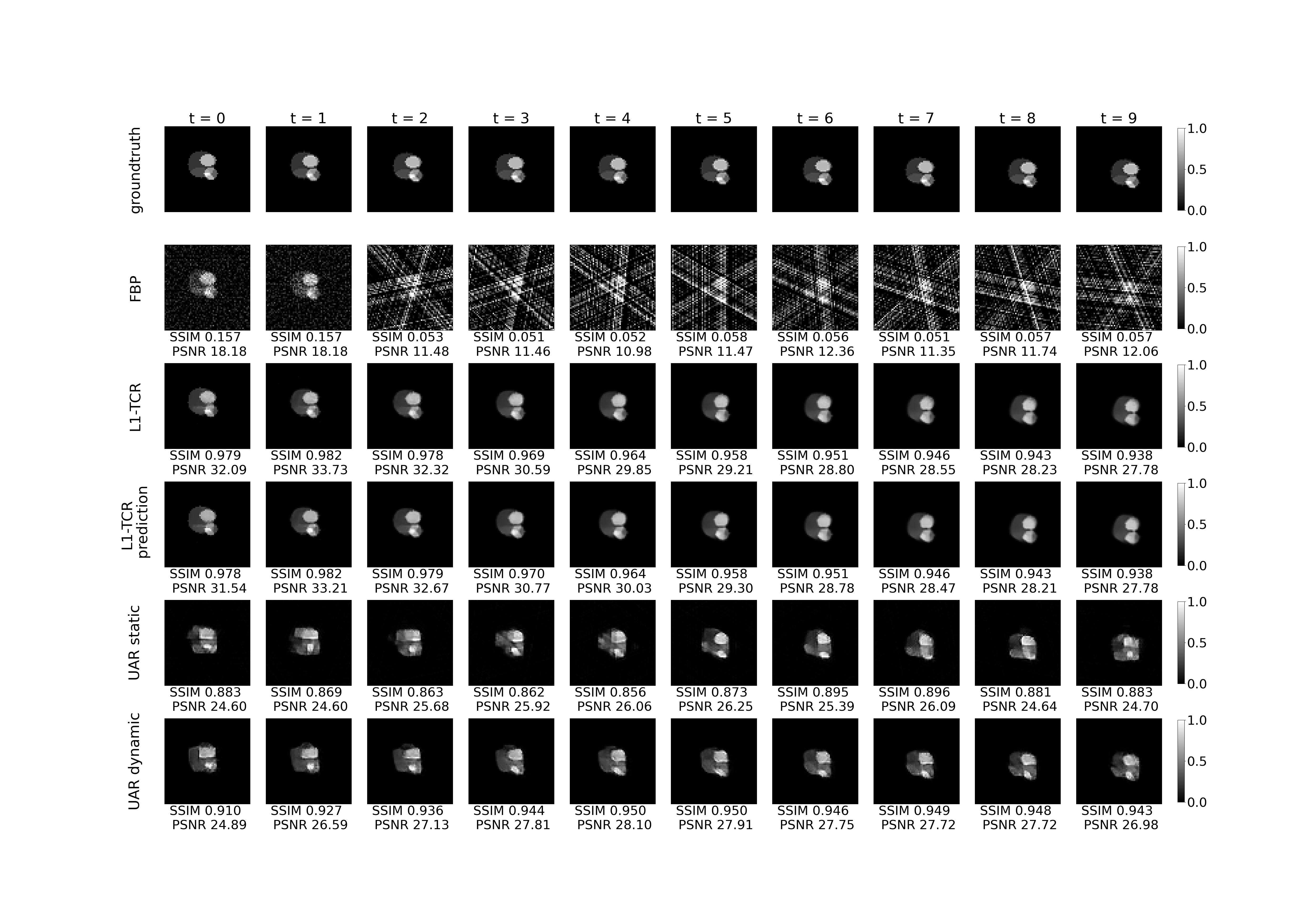}

}
    \caption{Results for a test phantom with 20 measurement angles for the initial time steps 0 and 1 for filtered backprojection (FBP), the $L^1$ transformer-based reconstruction (causality reconstruction) and corresponding prediction (causality prediction). For the remaining time steps measurements from 3 equidistantly rotating angles are taken. For the UAR reconstructions trained on static 2D phantoms (UAR static) and on dynamic 3D phantoms (UAR dynamic) measurements from 3 angles are taken for all time steps.}
    \label{fig:Phantom4}
\end{figure*}

In the following, we present results for initial reconstructions obtained with the initial baseline reconstruction network (before the adversarial training) 
and final reconstruction obtained with the final reconstruction network (after the adversarial training), see Section~\ref{sec:TrainingDetailsUAR}. 
The performance on the test data set can be seen in Table~\ref{tab:UAR} and the reconstructions of one phantom of the test data set are shown in the bottom rows of Figure~\ref{fig:Phantom4}.

\begin{figure*}[htbp]
    \centering
 
        \includegraphics[width=\textwidth,  trim=150 150 150 200,
  clip]{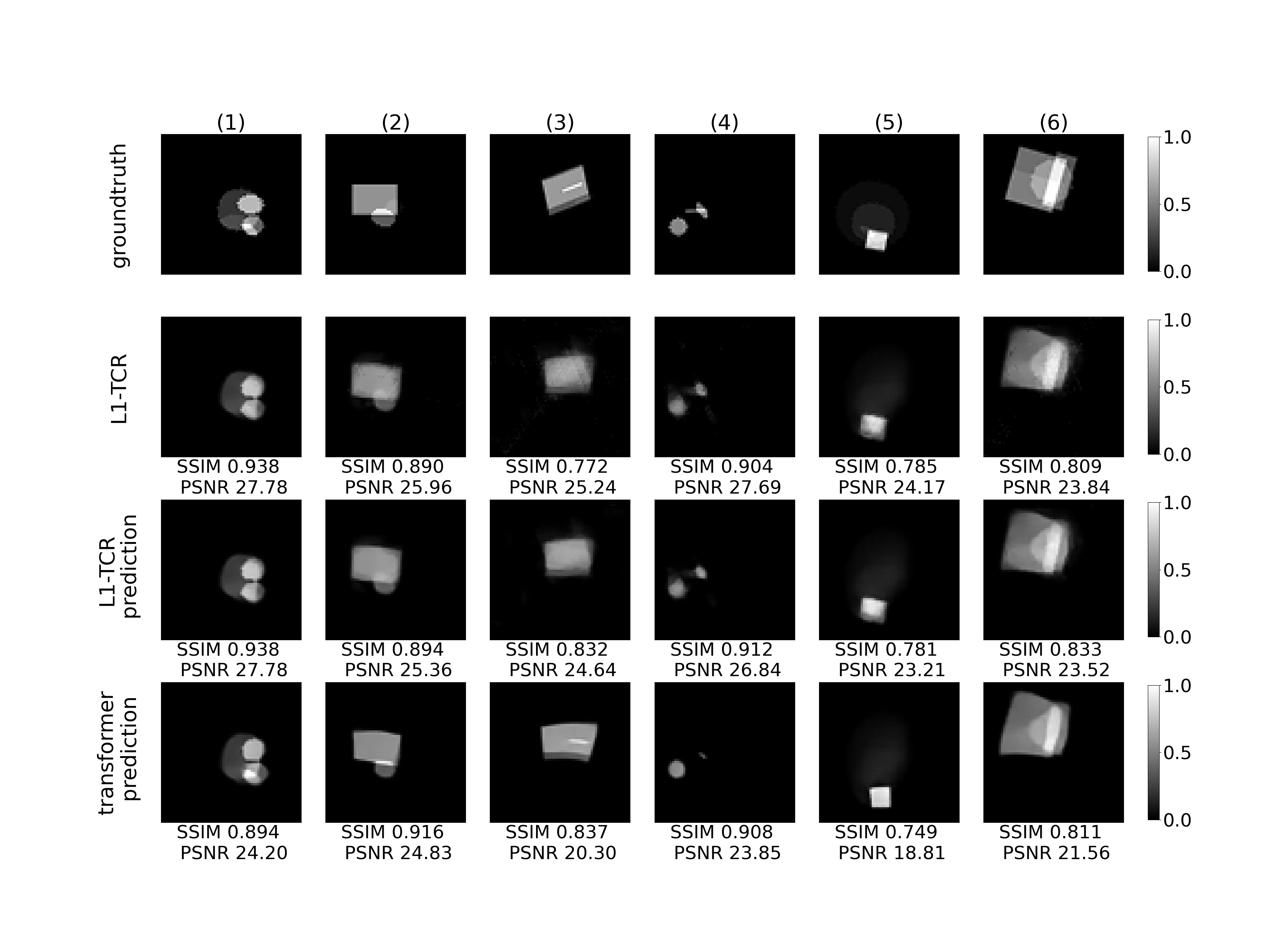}

    \caption{Results for last time step of several test phantoms with 3 measurement angles. From top to bottom: groundtruth, the $L^1$-transformer based reconstruction obtained with the learned causality function, the associated prediction and the prediction obtained auto-regressively by only using the transformer with two input frames (Landweber reconstructions obtained from 20 measurement angles). }
    \label{fig:LastTimestep}
\end{figure*}

\begin{table}[htb]
\centering
\begin{tabular}{|c|cc|}
\hline
\textbf{Static} & \multicolumn{1}{c}{PSNR} & SSIM  \\
\hline
3 angles: initial reco & 26.72, std 1.6466 & 0.887, std 0.045  \\
\multicolumn{1}{|r|}{final reco} & 26.27, std 1.561  &  0.834, std 0.0692 \\
\hline
10 angles: initial reco & 30.09, std 1.3488  & 0.9425, std 0.0262  \\
\multicolumn{1}{|r|}{final reco}  & 30.26, std 1.8532 & 0.9479, std 0.0337  \\
\hline
\hline
\textbf{Dynamic} & \multicolumn{1}{c}{PSNR} & SSIM \\
\hline
3 angles: initial reco & 24.1, std 1.205 & 0.8702, std 0.049  \\
\multicolumn{1}{|r|}{final reco} & 27.91, std 1.1298  & 0.9274, std 0.0243  \\
\hline
10 angles: initial reco & 26.95, std 0.9455 &  0.9065, std 0.0408 \\
\multicolumn{1}{|r|}{final reco} & 30.29, std 1.4917 & 0.9569, std 0.0189 \\
\hline
\end{tabular}
\vspace{0.5ex}
\caption{Mean PSNR and SSIM with standard deviation (std) of initial and final UAR reconstructions with 3 and 10 angles per time step for the static and dynamic setting.}
\label{tab:UAR}
\end{table}

\subsubsection{TCR}
We first test the auto-regressive performance of our transformer by using the first two frames of a phantom as input and then auto-regressively generating the remaining $8$ frames. We test this using Landweber reconstructions obtained from 20 angles or ground truth data. In Table~\ref{tab:transformerPerformance} the mean performance in SSIM and PSNR for this can be seen for the whole sequence and for the single last frame.

\begin{table}[htb]
    \centering
    \begin{tabular}{|c |cc|}
    \hline
        input data quality & PSNR  & SSIM \\
        \hline
        Landweber: all frames & 29.08, std 2.0932 &  0.9357, std 0.0275\\
        \multicolumn{1}{|r|}{last frame} & 24.71, std 2.9685 & 0.8875, std 0.0474 \\
        \hline
        ground truth, all frames & 29.07, std 2.0844 &  0.9518, std 0.0209 \\
        \multicolumn{1}{|r|}{last frame} & 25.44, std 3.354 & 0.9142, std 0.0417 \\
        \hline
    \end{tabular}
    \vspace{0.5ex}
    \caption{Mean PSNR and SSIM with standard deviation (std) of predicted sequences by the spatial-temporal transformer model with two input frames for the whole sequence and the last frame only.}
    \label{tab:transformerPerformance}
\end{table}

Next, we apply TCR to the test data set with $20$ angles for the first two time steps and $3$ or $10$ equidistantly spaced angles for the remaining steps.
 The mean PSNR and SSIM of the whole sequence and the last frame can be found in Table~\ref{tab:L1L2transformerPerformance}.

\begin{table}[htb]
    \centering
    \begin{tabular}{|c |cc|}
    \hline
        \textbf{$L^1$-TCR} & PSNR  & SSIM \\
        \hline
         3 angles, all frames & 29.63, std 1.9701 &  0.9296, std 0.0381\\
         \multicolumn{1}{|r|}{last frame} & 26.61, std 2.2742 & 0.8885, std 0.0565 \\
         \hline
       10 angles, all frames & 31.19, std 1.6231 &  0.943, std 0.0286 \\
       \multicolumn{1}{|r|}{last frame} & 29.92, std 1.6504 & 0.928, std 0.0339 \\
       \hline
       \hline
          \textbf{$L^2$-TCR}     & PSNR  & SSIM \\
        \hline
         3 angles, all frames & 28.76, std 1.7157&  0.7252, std 0.0708\\
        \multicolumn{1}{|r|}{last frame} & 26.17, std 2.2534&  0.6728, std 0.1136\\
         \hline
       10 angles, all frames & 29.45, std 1.3608 &  0.7501, std 0.0593 \\
       \multicolumn{1}{|r|}{last frame} & 28.04, std 1.4196&  0.7321, std 0.0621 \\
        \hline
    \end{tabular}
    \vspace{0.5ex}
    \caption{Mean PSNR and SSIM with standard deviation (std) of $L^1$- and $L^2$-TCR for the whole sequence and the last frame only.}
    \label{tab:L1L2transformerPerformance}
\end{table}

To compute the next prediction that serves as a prior in the minimization problem~\eqref{eq:DIP-NN-minProblem}, we use the previously computed reconstructions and not the output of the network, i.e., the auto-regressively generated predictions of the transformer, such that the predictions that we use for the minimization problem differ, see Figure~\ref{fig:LastTimestep}.

\subsection{Experimental data}
\begin{figure*}[htb]
\centering
\captionsetup[subfloat]{skip=0pt} 
\setlength{\tabcolsep}{0pt}

\begin{tabular}{c cccccccccc}
\vspace{-7ex}
\raisebox{-1.5\height}{\rotatebox{90}{\footnotesize reference}} &
\subfloat[]{%
\parbox[t]{0.095\textwidth}{\centering
{\footnotesize $t=0$}\par
\includegraphics[width=\linewidth,trim=100 0 100 0,clip]{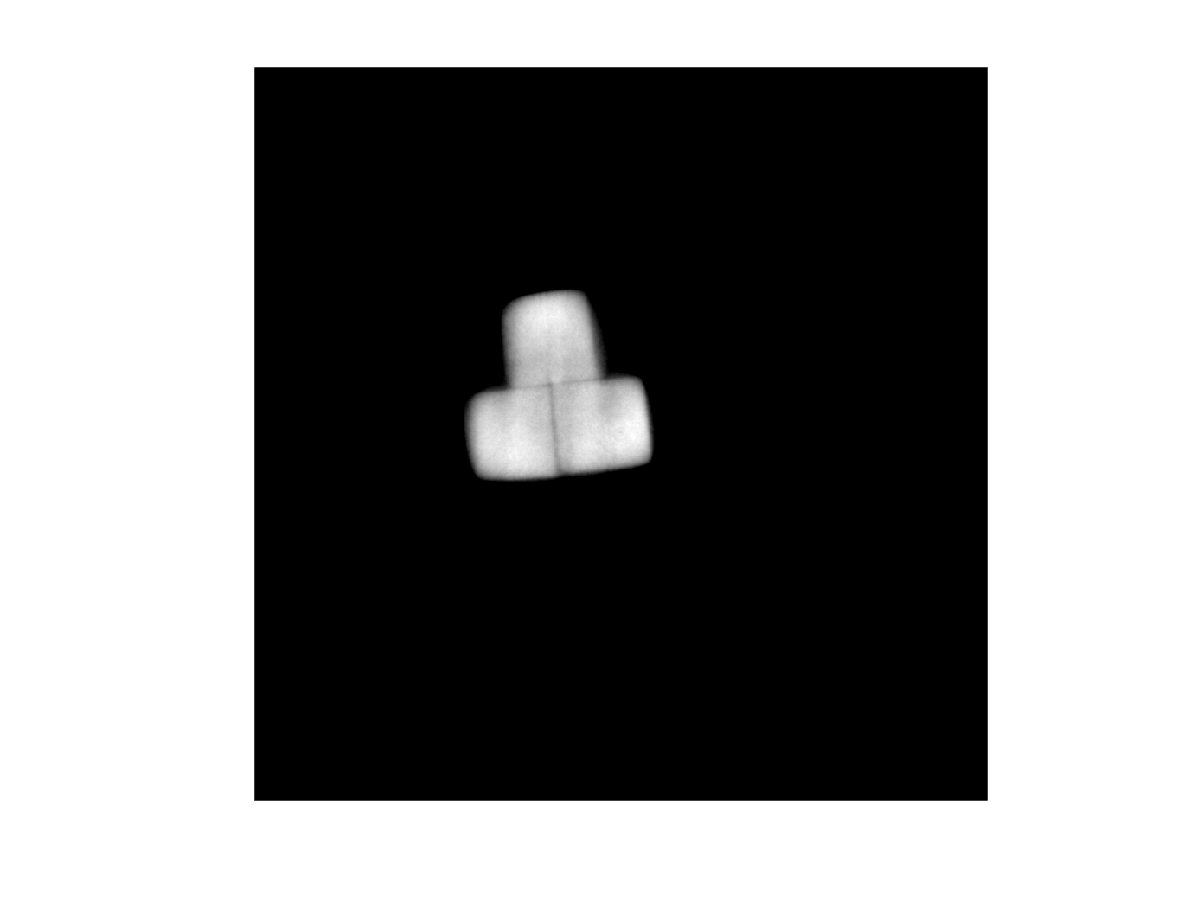}%
}%
} &
\subfloat[]{%
\parbox[t]{0.095\textwidth}{\centering
{\footnotesize $t=3$}\par
\includegraphics[width=\linewidth,trim=100 0 100 0,clip]{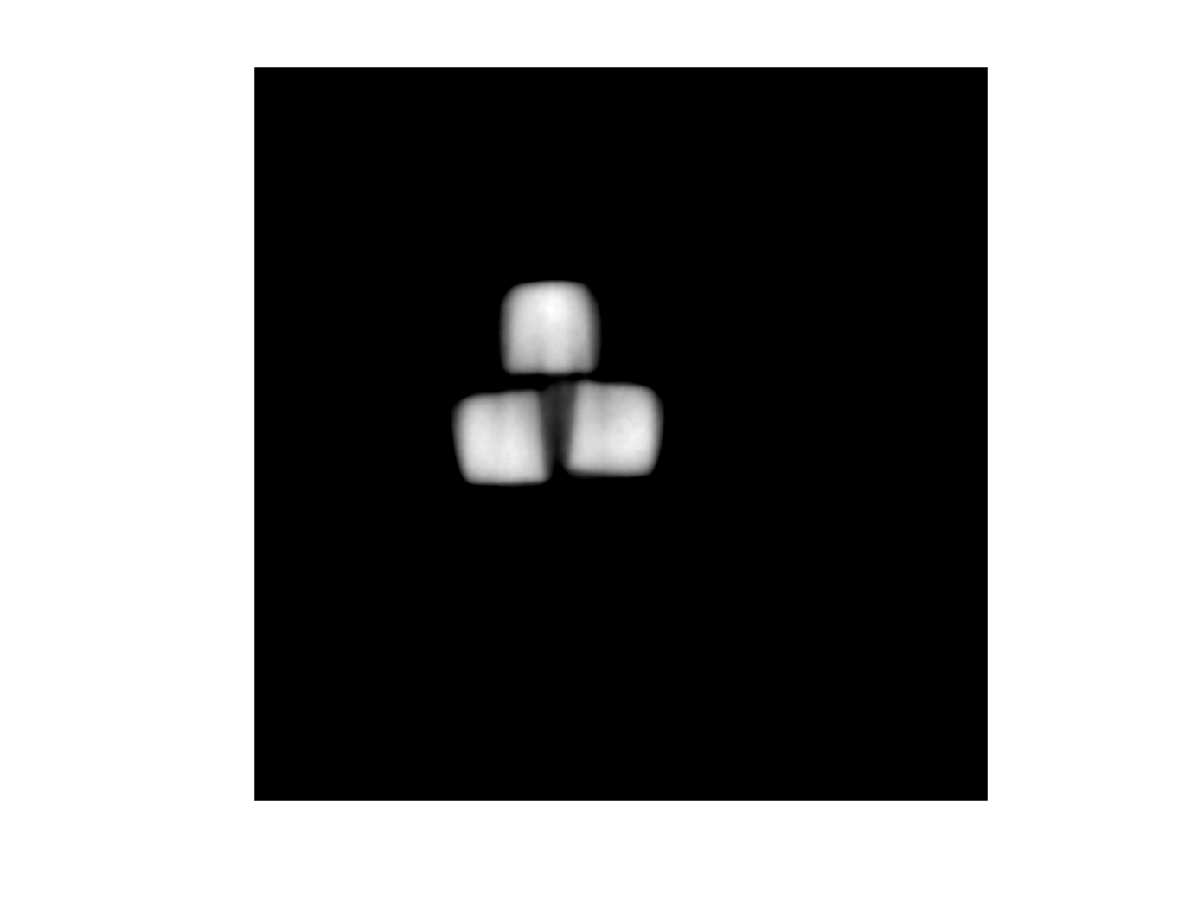}%
}%
} &
\subfloat[]{%
\parbox[t]{0.095\textwidth}{\centering
{\footnotesize $t=6$}\par
\includegraphics[width=\linewidth,trim=100 0 100 0,clip]{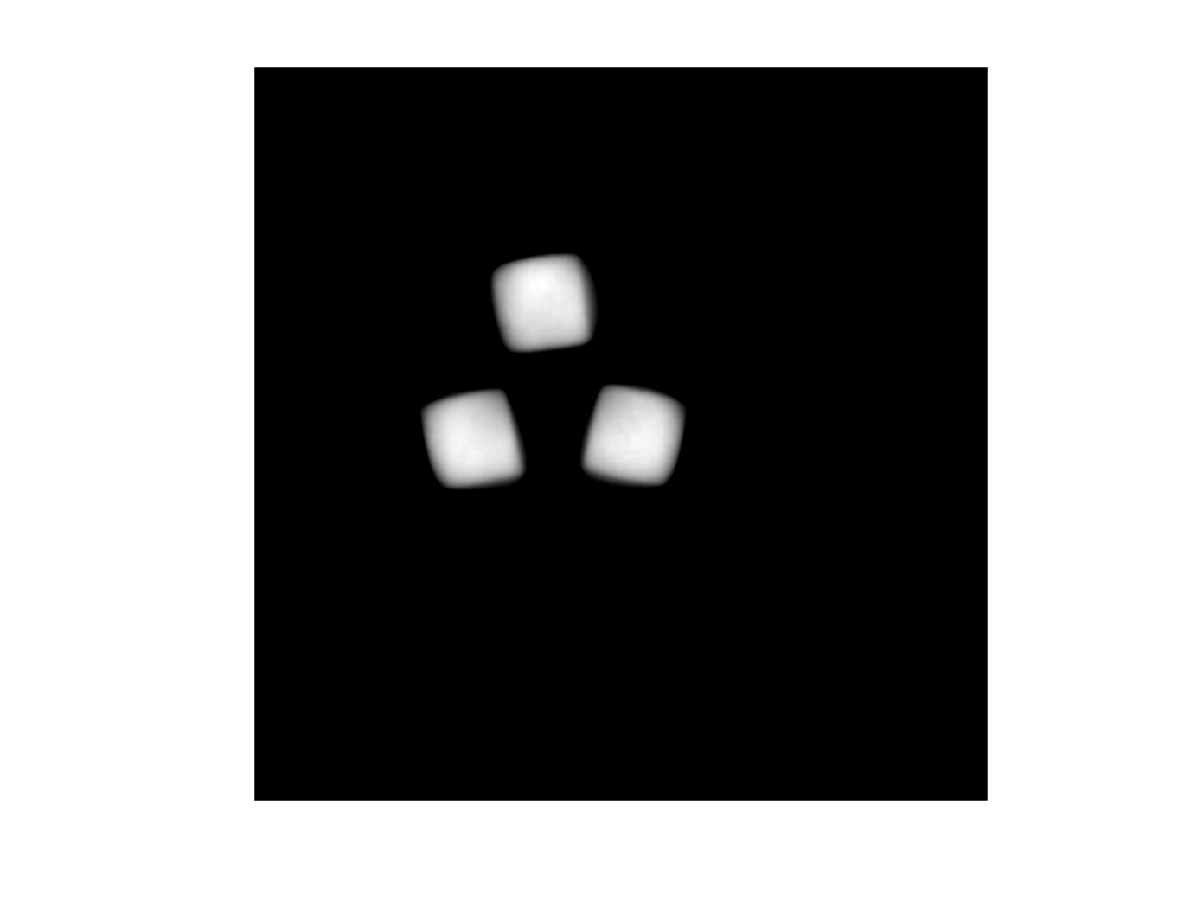}%
}%
} &
\subfloat[]{%
\parbox[t]{0.095\textwidth}{\centering
{\footnotesize $t=9$}\par
\includegraphics[width=\linewidth,trim=100 0 100 0,clip]{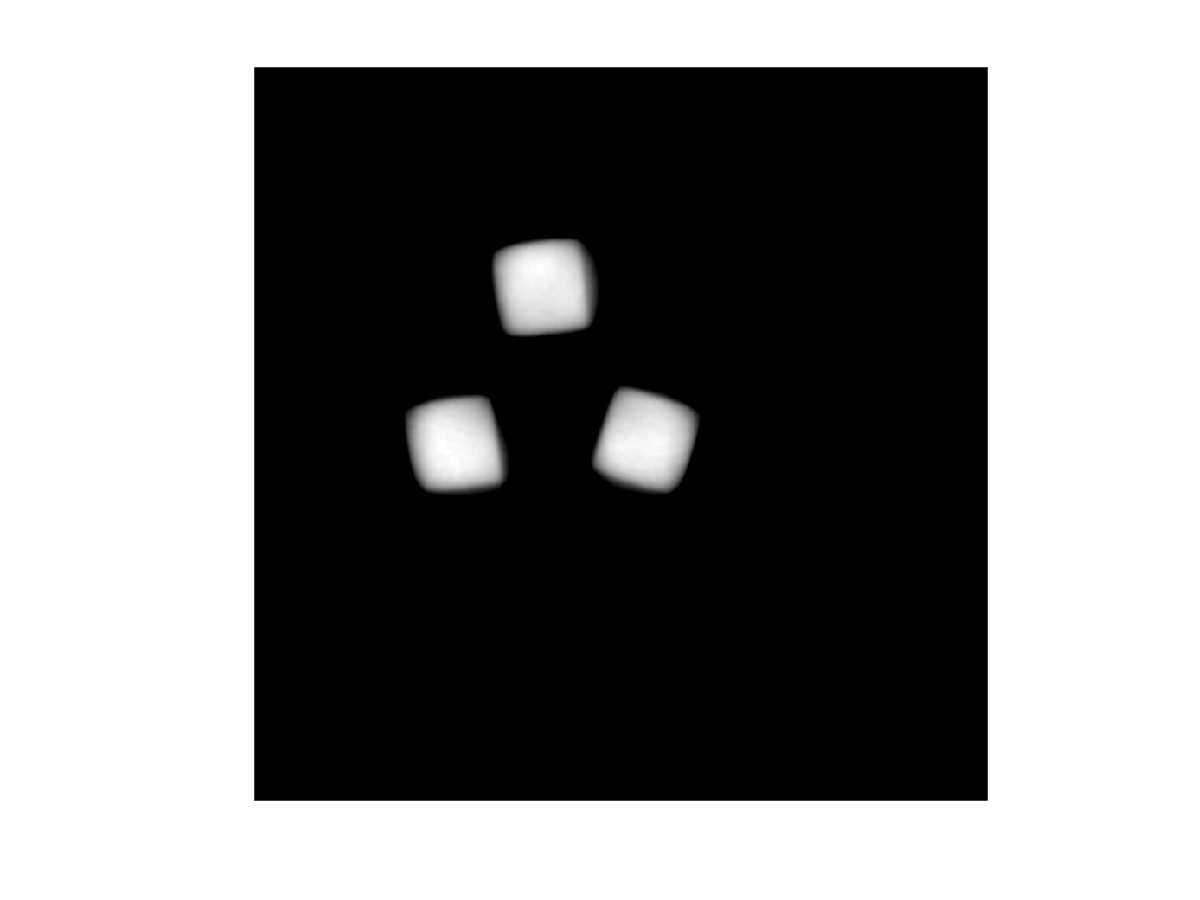}%
}%
}&
\subfloat[]{%
\parbox[t]{0.095\textwidth}{\centering
{\footnotesize $t=12$}\par
\includegraphics[width=\linewidth,trim=100 0 100 0,clip]{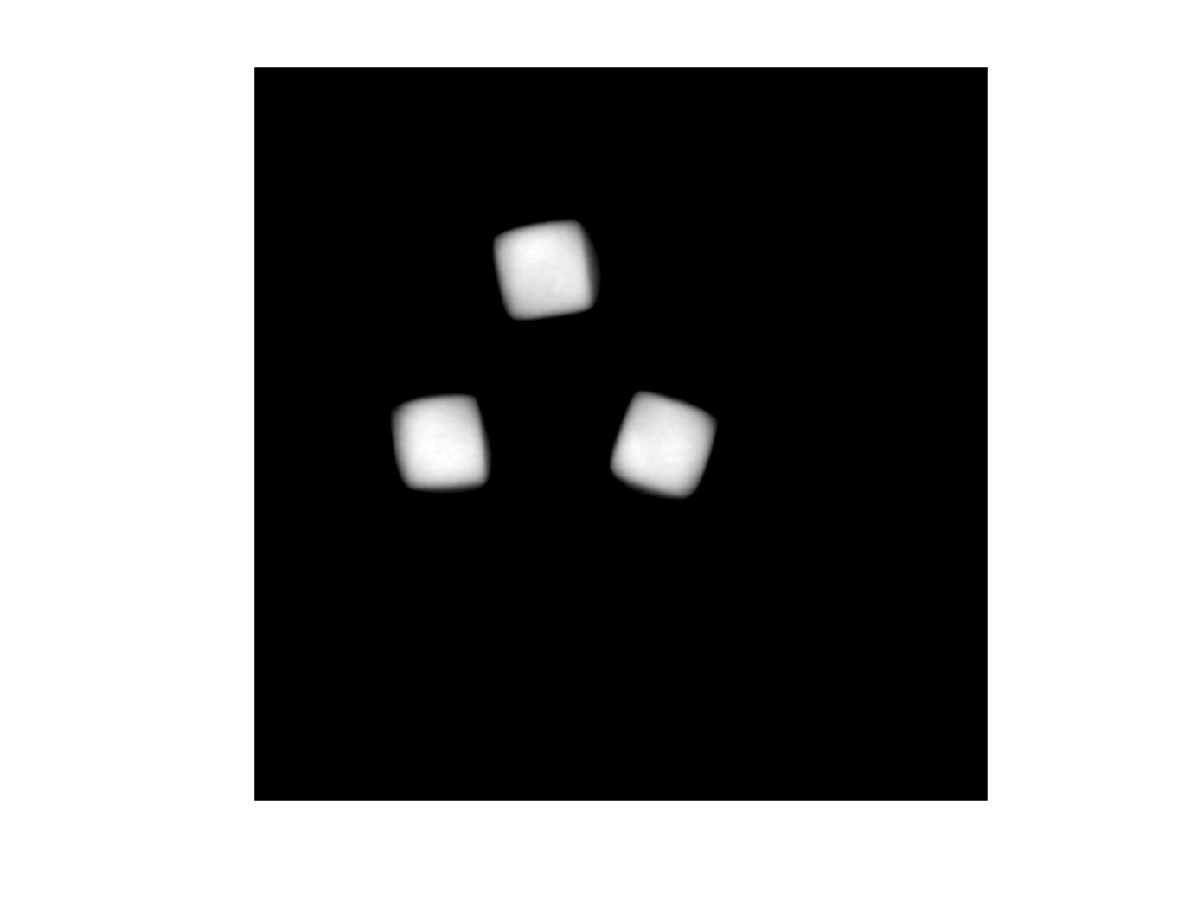}%
}%
} &
\subfloat[]{%
\parbox[t]{0.095\textwidth}{\centering
{\footnotesize $t=15$}\par
\includegraphics[width=\linewidth,trim=100 0 100 0,clip]{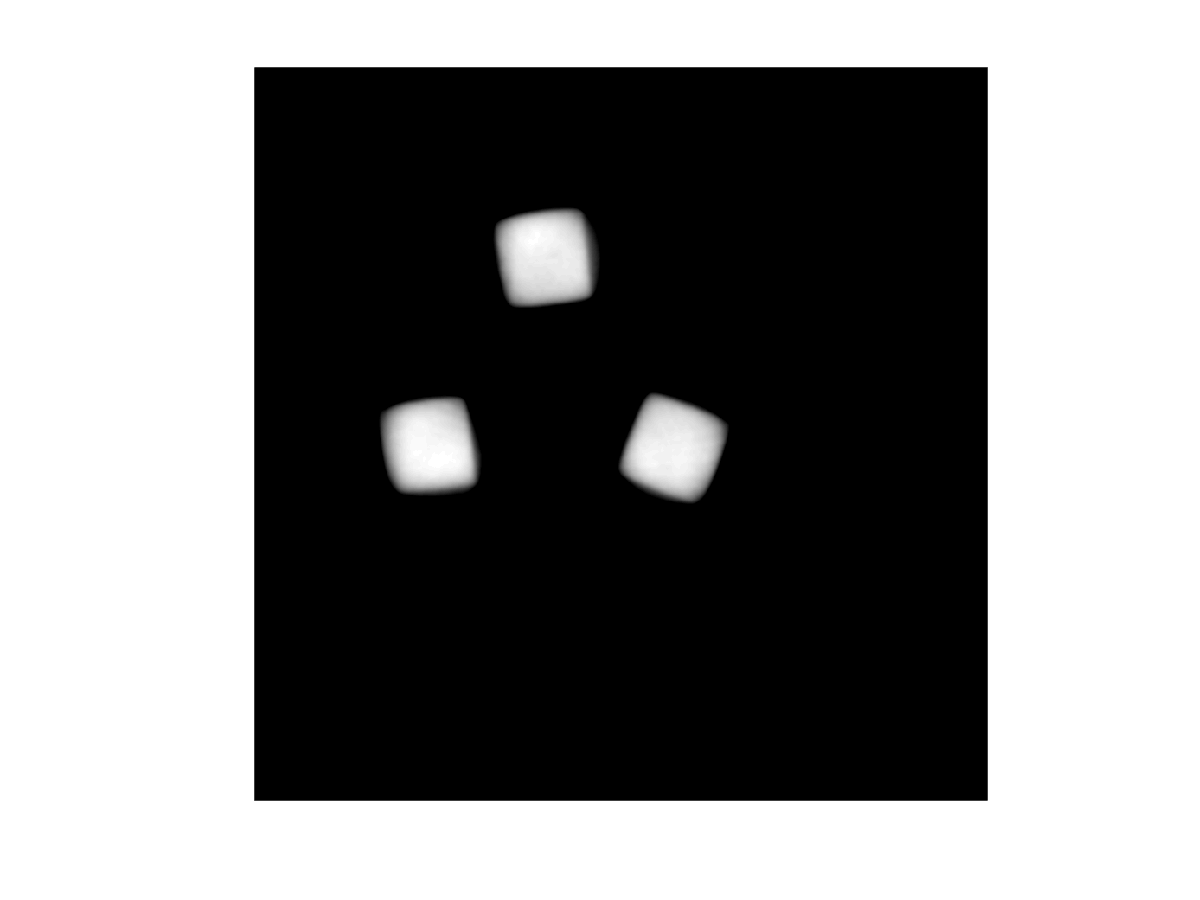}%
}%
} &
\subfloat[]{%
\parbox[t]{0.095\textwidth}{\centering
{\footnotesize $t=18$}\par
\includegraphics[width=\linewidth,trim=100 0 100 0,clip]{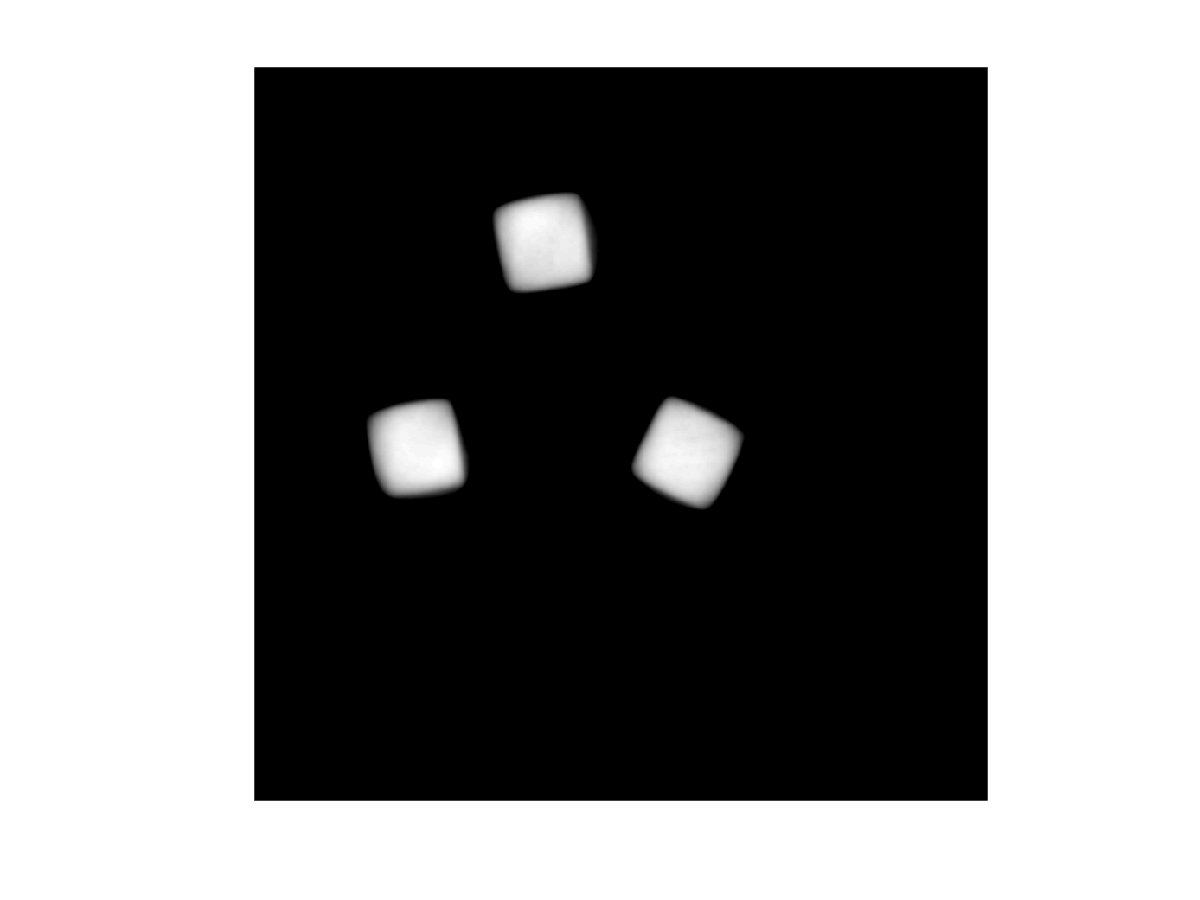}%
}%
} &
\subfloat[]{%
\parbox[t]{0.095\textwidth}{\centering
{\footnotesize $t=21$}\par
\includegraphics[width=\linewidth,trim=100 0 100 0,clip]{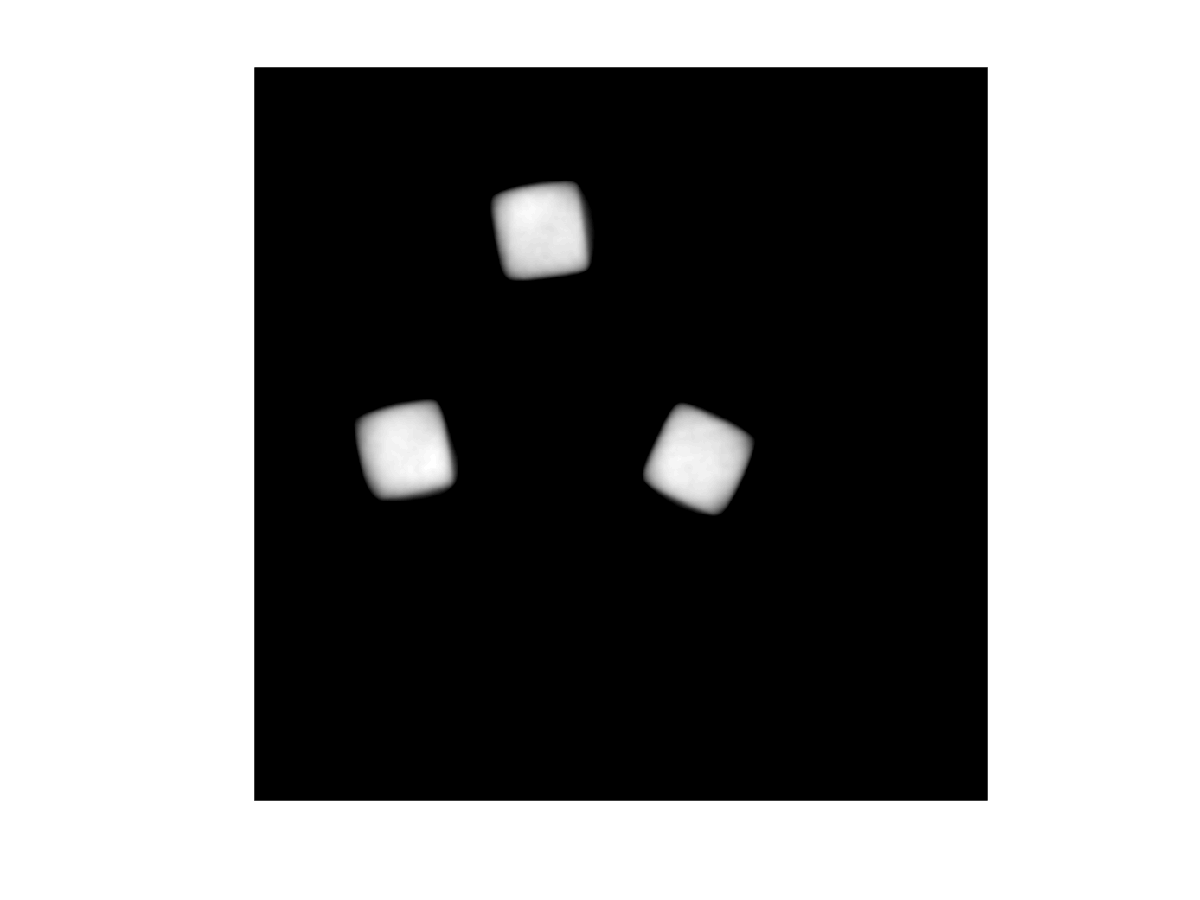}%
}%
} &
\subfloat[]{%
\parbox[t]{0.095\textwidth}{\centering
{\footnotesize $t=24$}\par
\includegraphics[width=\linewidth,trim=100 0 100 0,clip]{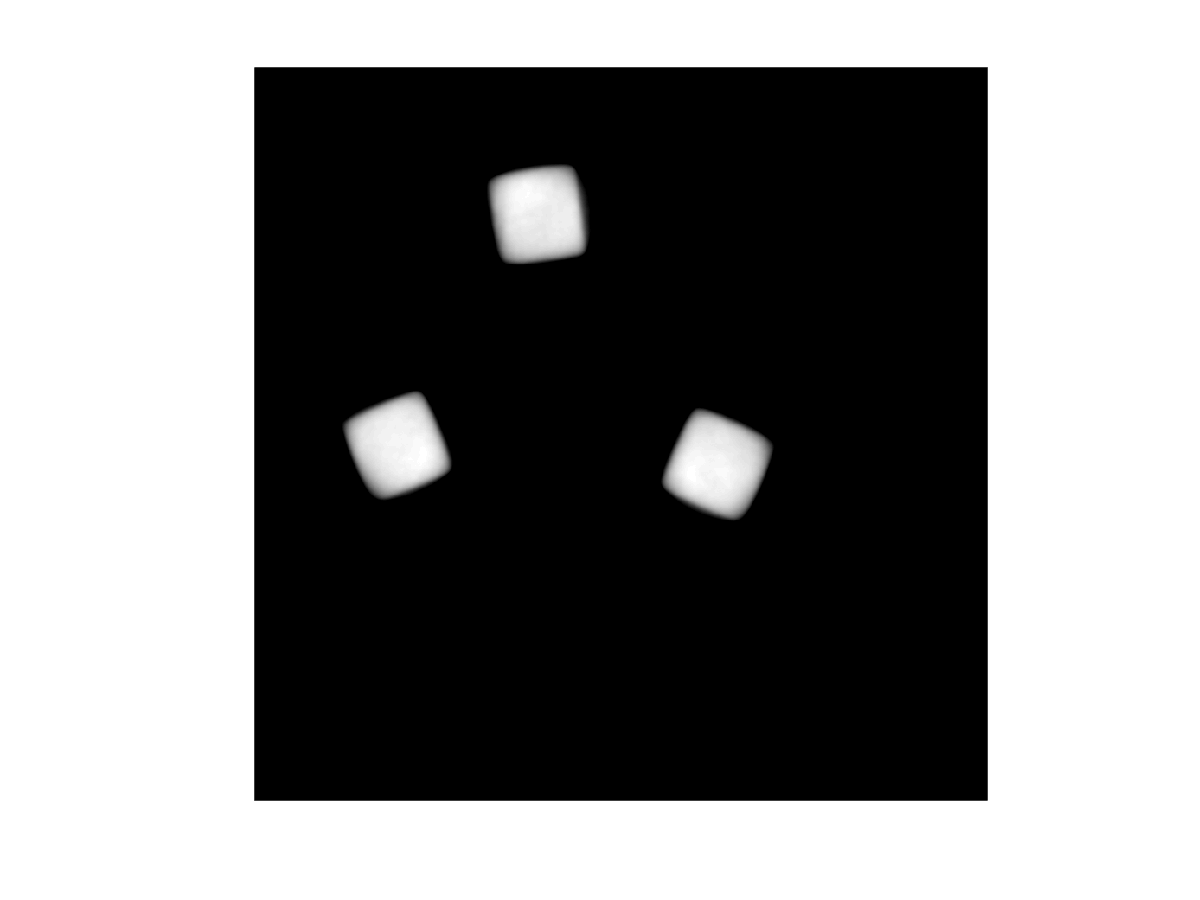}%
}%
} &
\subfloat[]{%
\parbox[t]{0.095\textwidth}{\centering
{\footnotesize $t=27$}\par
\includegraphics[width=\linewidth,trim=100 0 100 0,clip]{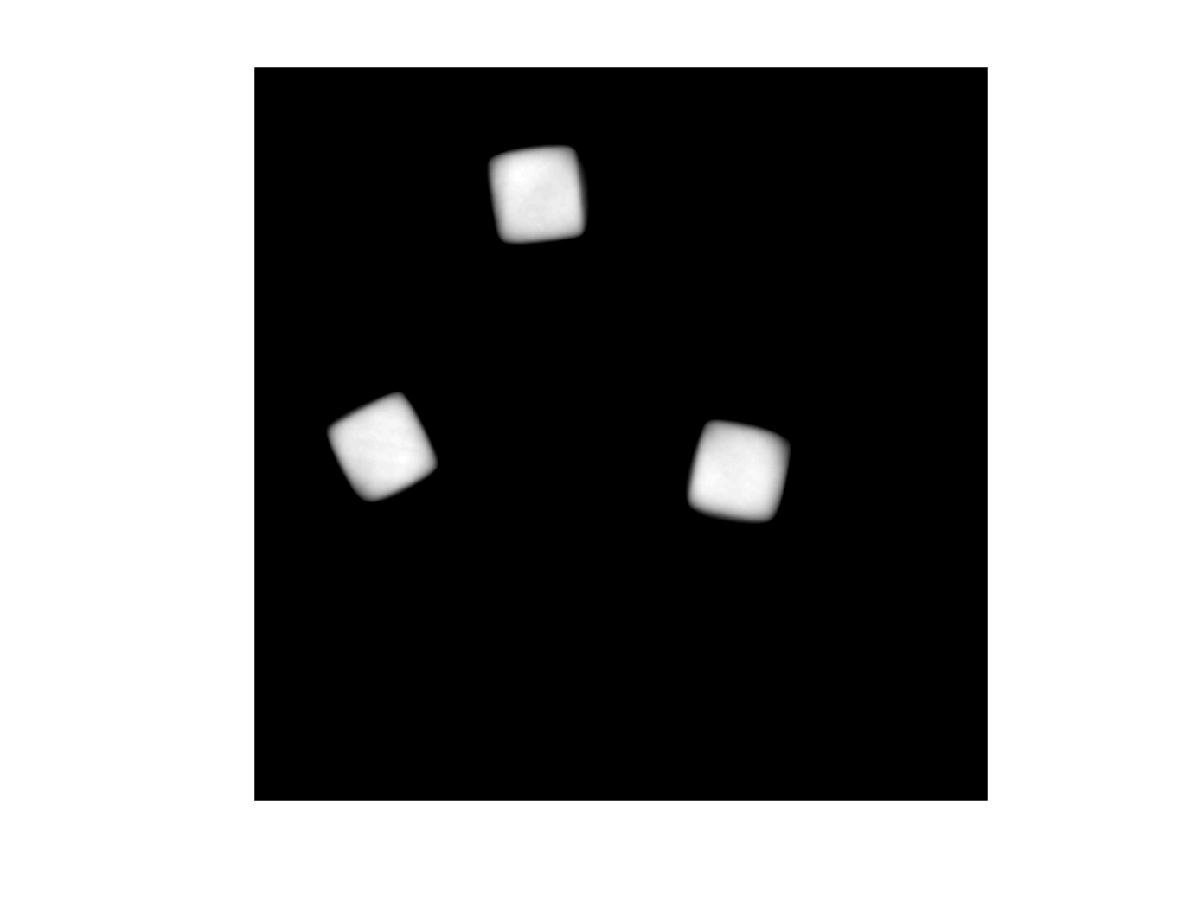}%
}%
} \\
\vspace{-7ex}
\rotatebox{90}{\, \scriptsize{$L^1$-TCR-pred} \quad \quad \scriptsize{$L^1$-TCR}} &
\multicolumn{10}{c}{%
  \subfloat[]{\centering
    \includegraphics[width=0.95\textwidth]{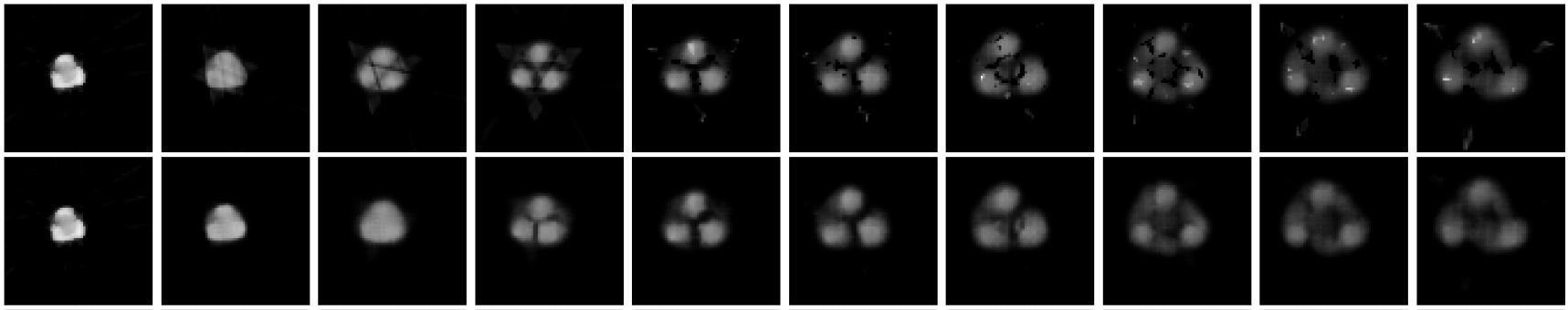}}%
} \\

\rotatebox{90}{\shortstack{\, \scriptsize{\, $L^1$-TV- \quad \quad \quad $L^1$-TV-} \\ \scriptsize{TCR-pred \qquad \quad \quad TCR}}} &
\multicolumn{10}{c}{%
  \subfloat[]{\centering
    \includegraphics[width=0.95\textwidth]{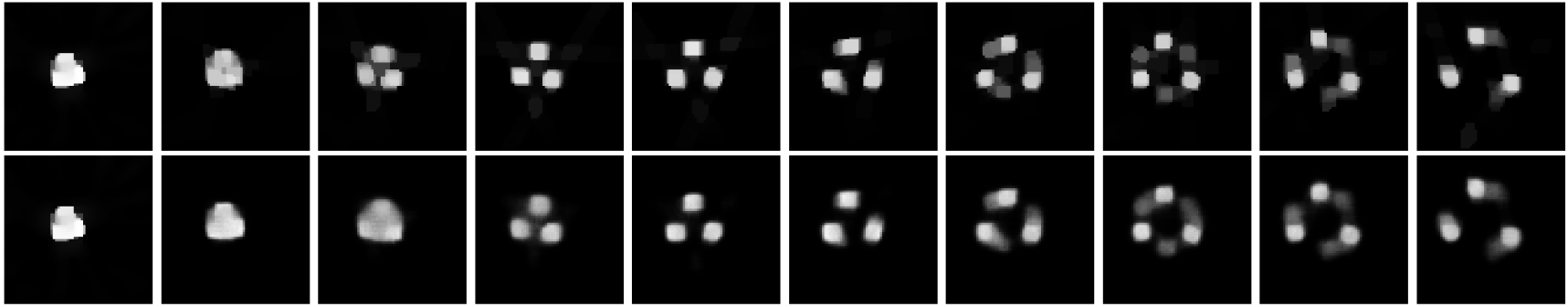}}%
} \\
\end{tabular}

\caption{From top to bottom: Reference solution for the rolling stones data set (reference), $L^1$ transformer based reconstruction for 10 initial angles and $3$ angles for the remaining time steps ($L^1$-reco),  corresponding transformer prediction ($L^1$-pred), $L^1$-TV transformer based reconstruction for 10 initial angles and $3$ angles for the remaining time steps ($L^1$-TV-reco), corresponding transformer prediction ($L^1$-TV-pred).}
\label{fig:rollingStones}
\end{figure*}

We apply TCR to the rolling stones data~\cite{Burger_2017}. Here, we use total variation (TV) regularization with a regularization parameter of $0.01$ instead of Landweber iteration to obtain the initial reconstructions for the first two time steps from $10$ measured angles each. These are then passed to the refinement model.
The results of $L^1$-TCR with an initial regularization parameter of $0.1$ for the first two time steps and a regularization parameter of $0.081$ for the remaining time steps are displayed in Figure~\ref{fig:rollingStones}. The predictions in Figure~\ref{fig:rollingStones} are again the predictions that the transformer produces with the previous reconstructions as input.

To improve the results, we add a TV-term additionally to the $L^1$-term, see~\eqref{eq:TVL1Min}. The initial regularization parameters are then $\alpha^{\text{init}}_{L^1} = 0.05$ and $\alpha^{\text{init}}_{TV} = 0.08$. For the remaining time steps, we choose $\alpha_{L^1} = 0.015$ and $\alpha_{TV} = 0.015$, see Figure~\ref{fig:rollingStones}.

\section{Discussion}
\label{sec:discussion}
\subsection{Numerical simulations}

\subsubsection{UAR}
It is striking that the static initial generator performs better than the dynamic one, whereas this changes for the final generator, see Table~\ref{tab:UAR}. Since the static data set is much smaller, it could be easier for the network to learn since it does not consider any time-dependent information and thus the adversarial training might not lead to a significant improvement. In contrast, in the dynamic case, the data set is bigger and contains dynamic information. Here, the adversarial training leads to a substantial improvement both for 3 and 10 angles.

\begin{figure}[htbp]
    \centering
        \includegraphics[width=0.5\textwidth,  trim=450 200 10 250,
  clip]{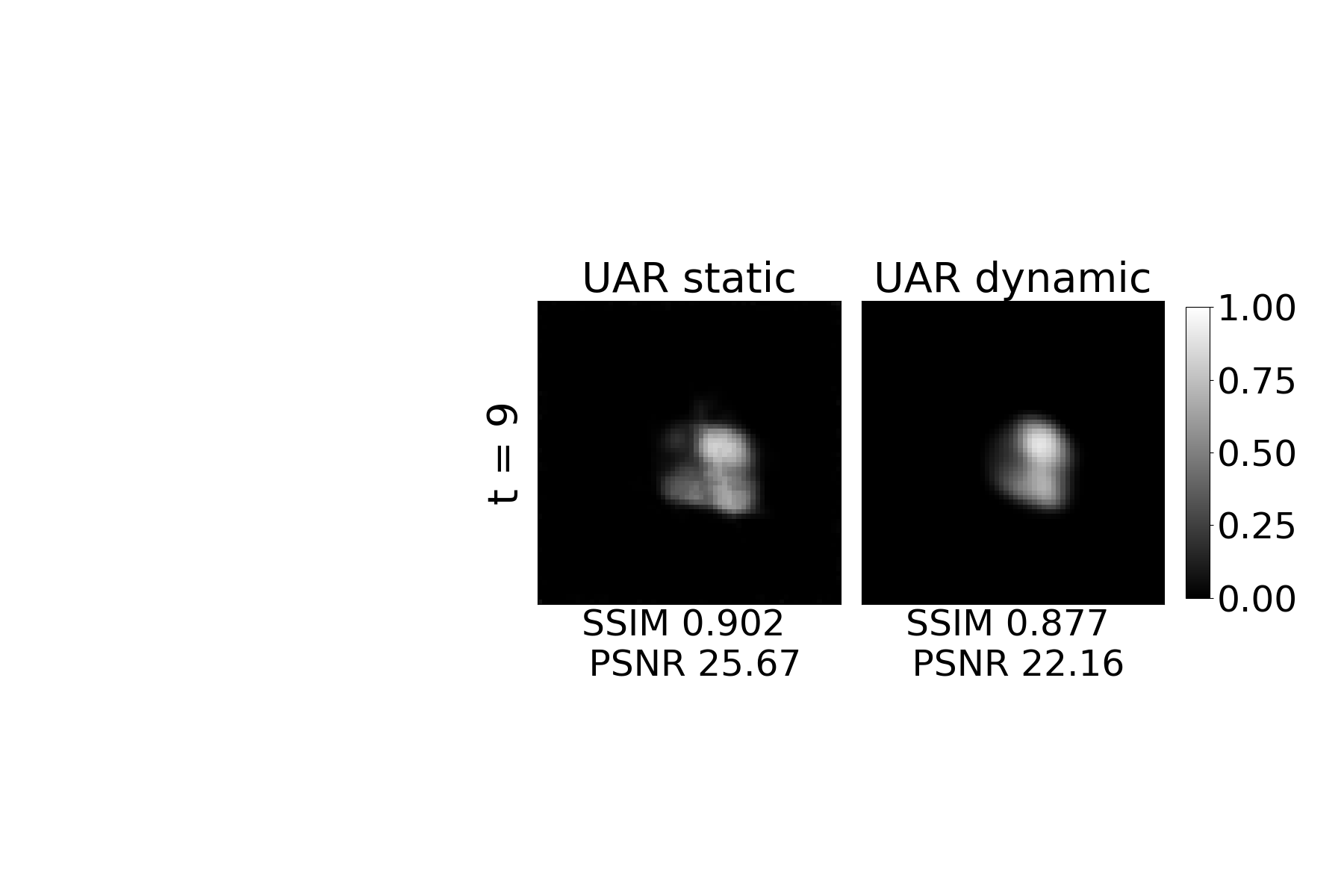}
    \caption{Last time step of the test phantom in Figure~\ref{fig:Phantom4} of initial baseline generator of UAR optimized using~\eqref{eq:lossGenDiscr} in the static (trained on 2D phantoms) and the dynamic case (trained on 3D phantoms).}
    \label{fig:initialGen}
\end{figure}

In addition, it is interesting that in the static UAR setting with 3 angles, the PSNR and SSIM of the initial generator are a bit higher than the ones of the final generator. However, if we look at the reconstructions, we see that the final reconstructions are much sharper than the initial ones and seem to be of better quality, see Figure~\ref{fig:initialGen}. Since $3$ angles per time step mean severe undersampling, and the model has no access to data from other time steps in the static setting, it is not surprising that the results behave unpredictable.

\subsubsection{TCR}

 The mean performance in SSIM and PSNR for the auto-regressively generated predictions from the spatial-temporal transformer, see Table~\ref{tab:transformerPerformance}, gives already reasonable results taking into account that we only have data for the first two time steps.

The $L^2$-TCR reconstructions obtained with Algorithm~\ref{algo_L2-TCR} can only improve the PSNR in the 10 angle case but not the SSIM compared to only using the transformer's predictions, see Table~\ref{tab:L1L2transformerPerformance}. $L^1$-TCR on the other hand improves the PSNR for both 3 and 10 angles, see Table~\ref{tab:L1L2transformerPerformance}. For 3 angles, the mean SSIM of all frames of the $L^1$-TCR reconstructions is a bit lower compared to the transformer's predictions. However, if we only consider the last frame, $L^1$-TCR improves the SSIM compared to the transformer's predictions, see Tables~\ref{tab:transformerPerformance} and~\ref{tab:L1L2transformerPerformance}. This shows that solving the variational problem for each time step adds data-consistency to the solution and reduces the error propagation of the transformer prediction model.

For 10 angles we obtain improved results regarding all measures. Already in the 3 angle case, we notice that the quality of the reconstructed last frame (in particular the respective PSNR) clearly improves compared to the transformer's predictions, see Figure~\ref{fig:LastTimestep}. This shows that the additional data for all time steps used in TCR helps to obtain better reconstructions up to the last time frame. The quality of the transformer predictions decreases over time, but by reconstructing with TCR this decrease in quality can be reduced. 

Note that in some cases the SSIM/PSNR of the predictions used for the solver are higher than the SSIM/PSNR of our final reconstructions. However, it still makes sense to use the reconstructions and not the auto-regressively generated predictions by the transformer since we only obtain these predictions based on the TCR reconstructions: When comparing SSIM/PSNR of TCR and the auto-regressively generated predictions, the reconstructions perform better in most cases, see Tables~\ref{tab:transformerPerformance} and~\ref{tab:L1L2transformerPerformance} and Figure~\ref{fig:LastTimestep}. This shows that TCR makes use of the additional data which we can also observe in a few examples.

\begin{figure}[htbp]
    \centering

        \includegraphics[width=0.5\textwidth,  trim=400 0 0 50,
  clip]{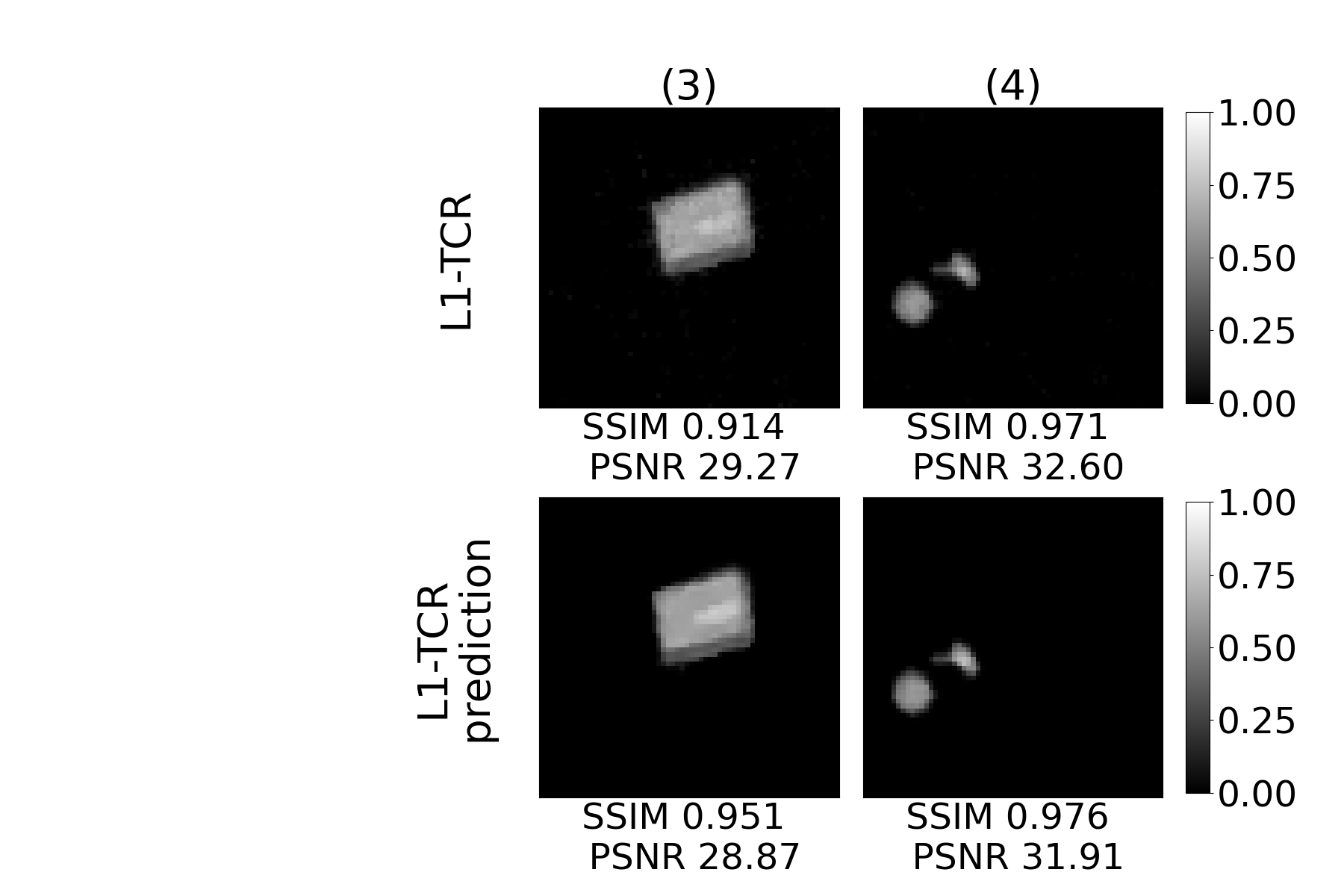}
    \caption{$L^1$ transformer based reconstruction and corresponding prediction for phantoms (3) and (4), see Figure~\ref{fig:LastTimestep}, with 10 measurement angles.}
    \label{fig:10angles}
\end{figure}

In some cases, we observe that features are not recovered by the transformer-based prediction only, but appear in the TCR reconstructions for the 3 angle case, see Phantoms (1), (4) and (5) in Figure~\ref{fig:LastTimestep}. For Phantom (1) the lower small circle is deformed in the transformer prediction, but this circle is clearly visible in the TCR reconstruction. The small line in Phantom (4) disappears in the transformer prediction. In the TCR reconstructions with 3 angles it is barely visible but does not disappear, whereas in the TCR reconstruction with 10 angles it is clearly visible, see Figure~\ref{fig:10angles}.
The rotation of the object in Phantom (5) cannot be seen in the transformer prediction but appears in the TCR reconstruction.

The transformer prediction for Phantom (6) is contorted. The TCR reconstruction appears a bit blurry, but the overall shape is reconstructed better. For Phantom (2) we observe a similar behavior. However, here the SSIM of the TCR reconstruction is worse than of the transformer prediction. This is not true for the PSNR and we see that the shape of the rectangle in the reconstruction is not as deformed as in the transformer case.

Phantom (3) is an example where the transformer's prediction does not work well. Also, the TCR reconstruction based on 3 angles does not recover the features properly. However, with 10 angles, the TCR reconstruction is improved, see Figure~\ref{fig:10angles}. This shows that additional data improves the final TCR reconstructions.

\subsubsection{Comparison of TCR and UAR}
As can be seen in Tables~\ref{tab:UAR} and~\ref{tab:L1L2transformerPerformance} the mean performance of TCR, our proposed variational approach using the learned causality function, outperforms the dynamic UAR approach in the 3 angles case in both PSNR and SSIM. In the 10 angle case, we can only improve the PSNR. However, the reconstructed time series seems more continuous in the TCR reconstructions compared to the UAR reconstructions even in cases where the SSIM of the UAR reconstructions is higher for some time steps, see Figure~\ref{fig:Phantom4}.

The quality of the reconstructions decreases over time in TCR, whereas this is not true for UAR since the UAR approach does not incorporate the causality principle and data for all time steps can be accessed at every time step.

TCR is flexible because it is independent of the number of time-points. The dynamic UAR approach could also be applied to a varying amount of time steps, but the model must be trained individually for each setting. The static setting can also be applied to a varying number of time points, however, the underlying dynamics are not integrated in this case.

\subsection{Experimental data}
The 'rolling stones data' undergoes a motion that does not belong to the possible training motions: In the training data set, the object as a whole undergoes an affine linear motion, whereas in the 'rolling stones data set', each of the three stones moves according to an individual - approximately affine linear - motion model. Furthermore, the transformer prediction model was only trained for up to 10 time points, while the 'rolling stones data' consist of 30 time points. Considering this, the reconstruction results for both $L^1$- and $L^1$-TV-regularization show that our trained model is able to generalize even when trained on a small dataset.

\section{Conclusions}
In this article, we extend the research on analytical solution methods for dynamic inverse problems in Lebesgue-Bochner-spaces~\cite{Burger_2024, Sarnighausen_2026}, which allow treating the spatial and temporal variable differently, to data-driven approaches.

We develop a hybrid regularization method called transformer causality regularization (TCR) that integrates a learned causality function as a prior into a classical variational regularization scheme. By combining the inductive bias of CNNs for imaging problems with the causal mask of transformers, our spatial-temporal transformer architecture can be seen as a transfer of the Bochner space setting to deep learning.
Since TCR reconstructions are obtained using classical variational regularization, the known solution methods and convergence properties can be directly transferred to our methods.

We compare this transformer-based approach to unrolled adversarial regularization (UAR)~\cite{mukherjee2021end}, a data-driven end-to-end method for solving linear inverse problems, by considering a static (trained on 2D phantoms) as well as a dynamic (trained on 3D phantoms) version of UAR.

$L^1$-TCR is able to outperform both the static and dynamic versions of UAR on the test data set. Furthermore, we apply TCR to experimentally measured data, which shows that the spatial-temporal-transformer can generalize to out-of-distribution data.

In comparison to the purely data-driven auto-regressively generated prediction of the spatial-temporal-transformer, we observe that the regularized TCR reconstructions can recover some features that are otherwise lost, motivating our interest in regularized solutions with theoretical convergence and regularization guaranties. In particular, by solving the variational problem we ensure data consistency of the regularized solutions.

This paper serves as a proof-of-concept for our developed method, and in future research this approach could be applied to more sophisticated indirect data, for example in medical imaging of moving objects or in nano-CT imaging, but also in other applications not connected to tomography. 

A drawback of our method is that we need initial reconstructions for the first two time steps. To reduce the amount of required data, the refinement model could be further improved.
On the one hand different types of reconstructions such as TV, $L^1$ or $L^2$ could be used during the training process to make the model more adaptable to different kinds of inputs. On the other hand, the quality of the reconstruction for the second time step could be gradually reduced in the training, i.e., such that the number of measurement angles for the second time step slowly decreases. This would lead to a more realistic setting, where we have several measurements available for the first time step but only a few measurements for all subsequent time steps.

In future research, it will be interesting to investigate the impact of different possibly nonlinear forward operators from other applied inverse problems such as magnetic resonance imaging or electrical impedance tomography. Since the transformer causality function is independent of the forward operator, no further training is necessary. However, when dealing with nonlinear problems, classical nonlinear solution methods must be applied instead of the linear methods described above.

\section*{Acknowledgements}
We would like to thank Carola-Bibiane Schönlieb and Constantin Pape for helpful discussions.

\bibliographystyle{IEEEtran}

\bibliography{source.bib}

@article{FISTA,
author = {Beck, Amir and Teboulle, Marc},
title = {A Fast Iterative Shrinkage-Thresholding Algorithm for Linear Inverse Problems},
journal = {SIAM Journal on Imaging Sciences},
volume = {2},
number = {1},
pages = {183-202},
year = {2009},
doi = {10.1137/080716542},

URL = { 
    
        https://doi.org/10.1137/080716542
    
    

},
eprint = { 
    
        https://doi.org/10.1137/080716542
    
    

}}

@INPROCEEDINGS{UNETR,

  author={Hatamizadeh, Ali and Tang, Yucheng and Nath, Vishwesh and Yang, Dong and Myronenko, Andriy and Landman, Bennett and Roth, Holger R. and Xu, Daguang},

  booktitle={2022 IEEE/CVF Winter Conference on Applications of Computer Vision (WACV)}, 

  title={UNETR: Transformers for 3D Medical Image Segmentation}, 

  year={2022},

  volume={},

  number={},

  pages={1748-1758},

  doi={10.1109/WACV51458.2022.00181}}

@inproceedings{17vaswaniAttention,
 author = {Vaswani, Ashish and Shazeer, Noam and Parmar, Niki and Uszkoreit, Jakob and Jones, Llion and Gomez, Aidan N and Kaiser, Lukasz and Polosukhin, Illia},
 booktitle = {Advances in Neural Information Processing Systems},
 editor = {I. Guyon and U. Von Luxburg and S. Bengio and H. Wallach and R. Fergus and S. Vishwanathan and R. Garnett},
 pages = {},
 publisher = {Curran Associates, Inc.},
 title = {Attention is All you Need},
 url = {https://proceedings.neurips.cc/paper_files/paper/2017/file/3f5ee243547dee91fbd053c1c4a845aa-Paper.pdf},
 volume = {30},
 year = {2017}
}

@article{24RoPE,
title = {RoFormer: Enhanced transformer with Rotary Position Embedding},
journal = {Neurocomputing},
volume = {568},
pages = {127063},
year = {2024},
issn = {0925-2312},
doi = {https://doi.org/10.1016/j.neucom.2023.127063},
url = {https://www.sciencedirect.com/science/article/pii/S0925231223011864},
author = {Jianlin Su and Murtadha Ahmed and Yu Lu and Shengfeng Pan and Wen Bo and Yunfeng Liu}
}

@inproceedings{20YunTransformerApproximator,
  author       = {Chulhee Yun and
                  Srinadh Bhojanapalli and
                  Ankit Singh Rawat and
                  Sashank J. Reddi and
                  Sanjiv Kumar},
  title        = {Are Transformers universal approximators of sequence-to-sequence functions?},
  booktitle    = {8th International Conference on Learning Representations, {ICLR} 2020,
                  Addis Ababa, Ethiopia, April 26-30, 2020},
  publisher    = {OpenReview.net},
  year         = {2020},
  url          = {https://openreview.net/forum?id=ByxRM0Ntvr},
  bibsource    = {dblp computer science bibliography, https://dblp.org}
}

@inproceedings{19adamw,
  author = {Ilya Loshchilov and Frank Hutter},
  title = {Decoupled Weight Decay Regularization},
  booktitle = {International Conference on Learning Representations},
  year = {2019},
  url = {https://openreview.net/forum?id=Bkg6RiCqY7}
}

@inproceedings{16SGDR_cosineAnnealing,
  author       = {Ilya Loshchilov and
                  Frank Hutter},
  title        = {{SGDR:} Stochastic Gradient Descent with Warm Restarts},
  booktitle    = {5th International Conference on Learning Representations, {ICLR} 2017,
                  Toulon, France, April 24-26, 2017, Conference Track Proceedings},
  year         = {2017},
  url          = {https://openreview.net/forum?id=Skq89Scxx},
  bibsource    = {dblp computer science bibliography, https://dblp.org}
}

@inproceedings{
furuya2025transformers,
title={Transformers are Universal In-context Learners},
author={Takashi Furuya and Maarten V. de Hoop and Gabriel Peyr{\'e}},
booktitle={The Thirteenth International Conference on Learning Representations},
year={2025},
url={https://openreview.net/forum?id=6S4WQD1LZR}
}

@inproceedings{Lunz2018,
    author = {Lunz, Sebastian and \"{O}ktem, Ozan and Sch\"{o}nlieb, Carola-Bibiane},
    title = {Adversarial regularizers in inverse problems},
    year = {2018},
    publisher = {Curran Associates Inc.},
    address = {Red Hook, NY, USA},
    booktitle = {Proceedings of the 32nd International Conference on Neural Information Processing Systems},
    pages = {8516–8525},
    numpages = {10},
    location = {Montr\'{e}al, Canada},
    series = {NIPS'18}
    }

@inproceedings{2021Mukherjee,
author = {Mukherjee, Subhadip and Carioni, Marcello and \"{O}ktem, Ozan and Sch\"{o}nlieb, Carola-Bibiane},
title = {End-to-end reconstruction meets data-driven regularization for inverse problems},
year = {2021},
isbn = {9781713845393},
publisher = {Curran Associates Inc.},
address = {Red Hook, NY, USA},
booktitle = {Proceedings of the 35th International Conference on Neural Information Processing Systems},
articleno = {1638},
numpages = {13},
series = {NIPS '21}
}

@Inbook{Bauschke2017Proximity,
author="Bauschke, Heinz H.
and Combettes, Patrick L.",
title="Proximity Operators",
bookTitle="Convex Analysis and Monotone Operator Theory in Hilbert Spaces",
year="2017",
publisher="Springer International Publishing",
address="Cham",
pages="413--446",
isbn="978-3-319-48311-5",
doi="10.1007/978-3-319-48311-5_24",
url="https://doi.org/10.1007/978-3-319-48311-5_24"
}

@article{14ParikhProximalAlgos,
author = {Parikh, Neal and Boyd, Stephen},
title = {Proximal Algorithms},
year = {2014},
issue_date = {Jan 2014},
publisher = {Now Publishers Inc.},
address = {Hanover, MA, USA},
volume = {1},
number = {3},
issn = {2167-3888},
url = {https://doi.org/10.1561/2400000003},
doi = {10.1561/2400000003},
journal = {Found. Trends Optim.},
month = jan,
pages = {127–239},
numpages = {116}
}

@article{2011ChambollPDHGAlgo,
author = {Chambolle, Antonin and Pock, Thomas},
title = {A First-Order Primal-Dual Algorithm for Convex Problems with Applications to Imaging},
year = {2011},
issue_date = {May       2011},
publisher = {Kluwer Academic Publishers},
address = {USA},
volume = {40},
number = {1},
issn = {0924-9907},
url = {https://doi.org/10.1007/s10851-010-0251-1},
doi = {10.1007/s10851-010-0251-1},
journal = {J. Math. Imaging Vis.},
month = may,
pages = {120–145},
numpages = {26}
}

@misc{jonas_adler_ODL,
  author       = {Jonas Adler and
                  Holger Kohr and
                  Ozan Öktem},
  title        = {Operator Discretization Library (ODL)},
  month        = jan,
  year         = {2017},
  doi          = {10.5281/zenodo.249479},
  url          = {https://doi.org/10.5281/zenodo.249479}
}

@article{Sarnighausen_2026,
doi = {10.1088/1361-6420/ae30f9},
url = {https://doi.org/10.1088/1361-6420/ae30f9},
year = {2026},
month = {jan},
publisher = {IOP Publishing},
volume = {42},
number = {1},
pages = {015008},
author = {Sarnighausen, Gesa and Hohage, Thorsten and Burger, Martin and Hauptmann, Andreas and Wald, Anne},
title = {Regularization for time-dependent inverse problems: geometry of {L}ebesgue–{B}ochner spaces and algorithms},
journal = {Inverse Problems}
}

@misc{albers2023timedependent,
      title={{Time-dependent parameter identification in a Fokker-Planck equation based magnetization model of large ensembles of nanoparticles}}, 
      author={Albers, H. and Kluth, T.},
      year={2023},
      eprint={2307.03560},
      howpublished={arXiv: 2307.03560},
      primaryClass={math.OC}
}

@article{bhw20,
doi = {10.1088/1361-6420/abb5e1},
url = {https://dx.doi.org/10.1088/1361-6420/abb5e1},
year = {2020},
publisher = {IOP Publishing},
volume = {36},
number = {12},
pages = {124001},
author = {Blanke, S. E. and Hahn, B. N. and Wald, A.},
title = {Inverse problems with inexact forward operator: iterative regularization and application in dynamic imaging},
journal = {Inverse Problems}
}

@article{Burger_2017,
doi = {10.1088/1361-6420/aa99cf},
url = {https://dx.doi.org/10.1088/1361-6420/aa99cf},
year = {2017},
publisher = {IOP Publishing},
volume = {33},
number = {12},
pages = {124008},
author = {Burger, M. and Dirks, H. and Frerking, L. and Hauptmann, A. and Helin, T. and Siltanen, S.},
title = {A variational reconstruction method for undersampled dynamic x-ray tomography based on physical motion models},
journal = {Inverse Problems}
}

@article{Burger_2024,
doi = {10.1088/1361-6420/ad5a35},
url = {https://dx.doi.org/10.1088/1361-6420/ad5a35},
year = {2024},
publisher = {IOP Publishing},
volume = {40},
number = {8},
pages = {085008},
author = {Burger, M. and Schuster, T. and Wald, A.},
title = {{Ill-posedness of time-dependent inverse problems in Lebesgue-Bochner spaces}},
journal = {Inverse Problems}
}

@book{Engl1996RegularizationOI,
  title={Regularization of Inverse Problems},
  author={H. W. Engl and M. Hanke and A. Neubauer},
  year={1996},
  url={https://api.semanticscholar.org/CorpusID:117875484},
  series = {{Mathematics and Its Applications}},
  publisher = {Springer Dordrecht}
}

@article{Goedeke_2023,
 author               = {G\"odeke, J. and Rigaud, G.},
 doi                  = {10.1088/1361-6420/acb2ed},
 journal              = {Inverse Problems},
 number               = {3},
 pages                = {034004},
 publisher            = {IOP Publishing},
 title                = {Imaging based on {C}ompton scattering: model uncertainty and data-driven reconstruction methods},
 url                  = {https://dx.doi.org/10.1088/1361-6420/acb2ed},
 volume               = {39},
 year                 = {2023},
 }

@article{Gris_2020,
 author               = {B. Gris and C. Chen and O. Öktem},
 doi                  = {10.1088/1361-6420/ab5832},
 journal              = {Inverse Problems},
 number               = {2},
 pages                = {025001},
 publisher            = {IOP Publishing},
 title                = {Image reconstruction through metamorphosis},
 url                  = {https://dx.doi.org/10.1088/1361-6420/ab5832},
 volume               = {36},
 year                 = {2020},
 }

@article{Hahn_2014,
 author               = {B. Hahn},
 doi                  = {10.1088/0266-5611/30/3/035008},
 journal              = {Inverse Problems},
 number               = {3},
 pages                = {035008},
 publisher            = {IOP Publishing},
 title                = {Efficient algorithms for linear dynamic inverse problems with known motion},
 url                  = {https://dx.doi.org/10.1088/0266-5611/30/3/035008},
 volume               = {30},
 year                 = {2014},
 }

@article{bh2014,
 author               = {Hahn, B.},
 journal              = {Journal of Inverse and Ill-posed Problems},
 number               = {3},
 pages                = {323-339},
 title                = {Reconstruction of dynamic objects with affine deformations in computerized tomography},
 volume               = {22},
 year                 = {2014},
 }

@article{bh2017,
 author               = {Hahn, B.},
 day                  = {04},
 doi                  = {10.1007/s11220-017-0159-6},
 issn                 = {1557-2072},
 journal              = {Sensing and Imaging},
 number               = {1},
 pages                = {10},
 title                = {Motion Estimation and Compensation Strategies in Dynamic Computerized Tomography},
 url                  = {https://doi.org/10.1007/s11220-017-0159-6},
 volume               = {18},
 year                 = {2017},
 }

@article{hashimoto2019dynamic,
 author               = {Hashimoto, F. and Ohba, H. and Ote, K. and Teramoto, A. and Tsukada, H.},
 journal              = {IEEE access},
 pages                = {96594--96603},
 publisher            = {IEEE},
 title                = {{Dynamic PET image denoising using deep convolutional neural networks without prior training datasets}},
 volume               = {7},
 year                 = {2019},
 }

@incollection{Hauptmann2021,
 address              = {Cham},
 author               = {Hauptmann, A. and {\"O}ktem, O. and Sch{\"o}nlieb, C.},
 booktitle            = {Handbook of Mathematical Models and Algorithms in Computer Vision and Imaging: Mathematical Imaging and Vision},
 doi                  = {10.1007/978-3-030-03009-4_83-1},
 editor               = {Chen, K. and Sch{\"o}nlieb, C.-B. and Tai, X.-C. and Younces, L.},
 isbn                 = {978-3-030-03009-4},
 pages                = {1--31},
 publisher            = {Springer International Publishing},
 title                = {Image Reconstruction in Dynamic Inverse Problems with Temporal Models},
 url                  = {https://doi.org/10.1007/978-3-030-03009-4_83-1},
 year                 = {2021},
 }

@article{kaltenbacher17,
 author               = {Kaltenbacher, B.},
 journal              = {Inverse Problems},
 number               = {6},
 pages                = {064002},
 title                = {All-at-once versus reduced iterative methods for time dependent inverse problems},
 url                  = {http://stacks.iop.org/0266-5611/33/i=6/a=064002},
 volume               = {33},
 year                 = {2017},
 }

@incollection{A1-KALTENBACHER;ET;AL:21,
 author               = {Kaltenbacher, B. and Nguyen, T. T. N. and Wald, A. and Schuster, T.},
 booktitle            = {Time-Dependent {{Problems}} in {{Imaging}} and {{Parameter Identification}}},
 doi                  = {10.1007/978-3-030-57784-1_13},
 editor               = {Kaltenbacher, B. and Schuster, T. and Wald, A.},
 isbn                 = {978-3-030-57784-1},
 langid               = {english},
 location             = {{Cham}},
 pages                = {377--412},
 publisher            = {{Springer International Publishing}},
 title                = {Parameter {{Identification}} for the {{Landau}}-{{Lifshitz}}-{{Gilbert Equation}} in {{Magnetic Particle Imaging}}},
 url                  = {https://doi.org/10.1007/978-3-030-57784-1_13},
 urldate              = {2022-07-01},
 year                 = {2021},
 }

@book{book_TDPIIP,
author="Kaltenbacher, B. and Schuster, T. and Wald, A.",
title="Time-dependent Problems in Imaging and Parameter Identification",
year="2021",
publisher="Springer International Publishing",
address="Cham",
isbn="978-3-030-57784-1",
doi="10.1007/978-3-030-57784-1",
url="https://doi.org/10.1007/978-3-030-57784-1"
}

@incollection{A1-KLEIN;ET;AL:21,
 author               = {Klein, R. and Schuster, T. and Wald, A.},
 booktitle            = {Time-Dependent {{Problems}} in {{Imaging}} and {{Parameter Identification}}},
 doi                  = {10.1007/978-3-030-57784-1_6},
 editor               = {Kaltenbacher, B. and Schuster, T. and Wald, A.},
 isbn                 = {978-3-030-57784-1},
 langid               = {english},
 location             = {{Cham}},
 pages                = {165--190},
 publisher            = {{Springer International Publishing}},
 title                = {Sequential {{Subspace Optimization}} for {{Recovering Stored Energy Functions}} in {{Hyperelastic Materials}} from {{Time-Dependent Data}}},
 url                  = {https://doi.org/10.1007/978-3-030-57784-1_6},
 urldate              = {2022-07-01},
 year                 = {2021},
 }

@book{C2-LAMBWAVES-BUCH:18,
 author               = {Lammering, R. and Grabbert, U. and Sinapius, M. and Schuster, T. and Wierach, P.},
 isbn                 = {978-3-319-49714-3},
 langid               = {english},
 publisher            = {{Springer}},
 title                = {{Lamb-Wave Based Structural Health Monitoring in Polymer Composites}},
 url                  = {https://openhsu.ub.hsu-hh.de/handle/10.24405/8678},
 urldate              = {2022-07-04},
 year                 = {2018},
 }

@article{Nitzsche22,
 author               = {Nitzsche, M. and Albers, H. and Kluth, T. and Hahn, B.},
 doi                  = {10.18416/IJMPI.2022.2203062},
 pages                = {165--190},
 publisher            = {Int J Mag Part Imag},
journal                = {International Journal on Magnetic Particle Imaging},
 title                = {Compensating model imperfections during image reconstruction via {R}esesop},
 url                  = {https://doi.org/10.18416/IJMPI.2022.2203062},
 volume               = {8},
 year                 = {2022},
 }

@book{Rieder:keineProblemeMitInversenProblemen,
author = {Rieder, A.},
year = {2003},
pages = {},
publisher="Vieweg+Teubner Verlag",
title = {Keine Probleme mit Inversen Problemen},
isbn = {978-3-528-03198-5},
doi = {10.1007/978-3-322-80234-7}
}

@article{Arridge_2019,
title={Solving inverse problems using data-driven models}, volume={28}, DOI={10.1017/S0962492919000059}, journal={Acta Numerica}, author={Arridge, Simon and Maass, Peter and Öktem, Ozan and Schönlieb, Carola-Bibiane}, year={2019}, pages={1–174}}

@Inbook{Haltmeier2021,
author="Haltmeier, Markus
and Nguyen, Linh",
editor="Chen, Ke
and Sch{\"o}nlieb, Carola-Bibiane
and Tai, Xue-Cheng
and Younces, Laurent",
title="Regularization of Inverse Problems by Neural Networks",
bookTitle="Handbook of Mathematical Models and Algorithms in Computer Vision and Imaging: Mathematical Imaging and Vision",
year="2021",
publisher="Springer International Publishing",
address="Cham",
pages="1--29",
isbn="978-3-030-03009-4",
doi="10.1007/978-3-030-03009-4_81-1",
url="https://doi.org/10.1007/978-3-030-03009-4_81-1"
}

@article{kobler2021total,
  title={Total deep variation: A stable regularization method for inverse problems},
  author={Kobler, Erich and Effland, Alexander and Kunisch, Karl and Pock, Thomas},
  journal={IEEE transactions on pattern analysis and machine intelligence},
  volume={44},
  number={12},
  pages={9163--9180},
  year={2021},
  publisher={IEEE}
}

@article{li2020nett,
  title={NETT: solving inverse problems with deep neural networks},
  author={Li, Housen and Schwab, Johannes and Antholzer, Stephan and Haltmeier, Markus},
  journal={Inverse Problems},
  volume={36},
  number={6},
  pages={065005},
  year={2020},
  publisher={IOP Publishing}
}

@article{mukherjee2021end,
  title={End-to-end reconstruction meets data-driven regularization for inverse problems},
  author={Mukherjee, Subhadip and Carioni, Marcello and {\"O}ktem, Ozan and Sch{\"o}nlieb, Carola-Bibiane},
  journal={Advances in Neural Information Processing Systems},
  volume={34},
  pages={21413--21425},
  year={2021}
}

@article{mukherjee2023learned,
  title={Learned reconstruction methods with convergence guarantees: A survey of concepts and applications},
  author={Mukherjee, Subhadip and Hauptmann, Andreas and {\"O}ktem, Ozan and Pereyra, Marcelo and Sch{\"o}nlieb, Carola-Bibiane},
  journal={IEEE Signal Processing Magazine},
  volume={40},
  number={1},
  pages={164--182},
  year={2023},
  publisher={IEEE}
}

@inproceedings{mukherjee2024data,
  title={Data-driven convex regularizers for inverse problems},
  author={Mukherjee, Subhadip and Dittmer, S{\"o}ren and Shumaylov, Zakhar and Lunz, Sebastian and {\"O}ktem, Ozan and Sch{\"o}nlieb, C-B},
  booktitle={ICASSP 2024-2024 IEEE International Conference on Acoustics, Speech and Signal Processing (ICASSP)},
  pages={13386--13390},
  year={2024},
  organization={IEEE}
}

@inproceedings{ahuja2023transformers,
  title={Transformers can learn to solve linear-inverse problems in-context},
  author={Ahuja, Kabir and Panwar, Madhur and Goyal, Navin},
  booktitle={NeurIPS 2023 Workshop on Deep Learning and Inverse Problems},
  year={2023}
}

@article{guo2022transformer,
  title={Transformer meets boundary value inverse problems},
  author={Guo, Ruchi and Cao, Shuhao and Chen, Long},
  journal={arXiv preprint arXiv:2209.14977},
  year={2022}
}

@ARTICLE{23ZhouNNformer,
  author={Zhou, Hong-Yu and Guo, Jiansen and Zhang, Yinghao and Han, Xiaoguang and Yu, Lequan and Wang, Liansheng and Yu, Yizhou},
  journal={IEEE Transactions on Image Processing}, 
  title={nnFormer: Volumetric Medical Image Segmentation via a 3D Transformer}, 
  year={2023},
  volume={32},
  number={},
  pages={4036-4045},
  doi={10.1109/TIP.2023.3293771}}

@InProceedings{21Petit,
author="Petit, Olivier
and Thome, Nicolas
and Rambour, Clement
and Themyr, Loic
and Collins, Toby
and Soler, Luc",
editor="Lian, Chunfeng
and Cao, Xiaohuan
and Rekik, Islem
and Xu, Xuanang
and Yan, Pingkun",
title="U-Net Transformer: Self and Cross Attention for Medical Image Segmentation",
booktitle="Machine Learning in Medical Imaging",
year="2021",
publisher="Springer International Publishing",
address="Cham",
pages="267--276",
isbn="978-3-030-87589-3"
}

@article{du2024inhomogeneous,
  title={Inhomogeneous media inverse scattering problem assisted by Swin transformer network},
  author={Du, Naike and Wang, Jing and Song, Rencheng and Xu, Kuiwen and Sun, Sheng and Ye, Xiuzhu},
  journal={IEEE Transactions on Microwave Theory and Techniques},
  volume={72},
  number={12},
  pages={6809--6820},
  year={2024},
  publisher={IEEE}
}

@inproceedings{muthukrishnan2023invrt,
  title={Invrt: Solving radar inverse problems with transformers},
  author={Muthukrishnan, Ramya and Goodwin, Justin and Kern, Adam and Vaska, Nathan and Caceres, Rajmonda S},
  booktitle={AAAI},
  volume={2},
  number={3},
  pages={4},
  year={2023}
}

@article{sun2025video,
  title={Video reconstruction through dynamic scattering media based on physics-informed spatio-temporal transformer},
  author={Sun, Peng and Wang, Canjin and Wang, Rijun and Pan, Shichao and Liu, Ji and Wang, Gao and Gong, Lukai and Gu, Lian},
  journal={Optics Express},
  volume={33},
  number={25},
  pages={52418--52432},
  year={2025},
  publisher={Optica Publishing Group}
}

@article{daubechies2004iterative,
  title={An iterative thresholding algorithm for linear inverse problems with a sparsity constraint},
  author={Daubechies, Ingrid and Defrise, Michel and De Mol, Christine},
  journal={Communications on Pure and Applied Mathematics: A Journal Issued by the Courant Institute of Mathematical Sciences},
  volume={57},
  number={11},
  pages={1413--1457},
  year={2004},
  publisher={Wiley Online Library}
}

@inproceedings{li2022mat,
  title={Mat: Mask-aware transformer for large hole image inpainting},
  author={Li, Wenbo and Lin, Zhe and Zhou, Kun and Qi, Lu and Wang, Yi and Jia, Jiaya},
  booktitle={Proceedings of the IEEE/CVF conference on computer vision and pattern recognition},
  pages={10758--10768},
  year={2022}
}

@ARTICLE{zhao22_transformerCNNDenoising,
  author={Zhao, Mo and Cao, Gang and Huang, Xianglin and Yang, Lifang},
  journal={IEEE Signal Processing Letters}, 
  title={Hybrid Transformer-CNN for Real Image Denoising}, 
  year={2022},
  volume={29},
  number={},
  pages={1252-1256},
  doi={10.1109/LSP.2022.3176486}}

@article{TIAN2024_CrossTransformerDenoising,
title = {A cross Transformer for image denoising},
journal = {Information Fusion},
volume = {102},
pages = {102043},
year = {2024},
issn = {1566-2535},
doi = {https://doi.org/10.1016/j.inffus.2023.102043},
url = {https://www.sciencedirect.com/science/article/pii/S1566253523003597},
author = {Chunwei Tian and Menghua Zheng and Wangmeng Zuo and Shichao Zhang and Yanning Zhang and Chia-Wen Lin}
}

@article{luetjen24,
  author={Lütjen, Tom and Schönfeld, Fabian and Oberacker, Alice and Leuschner, Johannes and Schmidt, Maximilian and Wald, Anne and Kluth, Tobias},
  journal={IEEE Transactions on Computational Imaging}, 
  title={Learning-Based Approaches for Reconstructions With Inexact Operators in nano{C}{T} Applications},
  year={2024},
  volume={10},
  number={},
  pages={522--534},
  doi={10.1109/TCI.2024.3380319}
}

@article{hauptmann2025convergent,
  title={Convergent regularization in inverse problems and linear plug-and-play denoisers},
  author={Hauptmann, Andreas and Mukherjee, Subhadip and Sch{\"o}nlieb, Carola-Bibiane and Sherry, Ferdia},
  journal={Foundations of Computational Mathematics},
  volume={25},
  number={4},
  pages={1087--1120},
  year={2025},
  publisher={Springer}
}

@article{tallman2020structural,
  title={Structural health and condition monitoring via electrical impedance tomography in self-sensing materials: a review},
  author={Tallman, Tyler N and Smyl, Danny J},
  journal={Smart Materials and Structures},
  volume={29},
  number={12},
  pages={123001},
  year={2020},
  publisher={IOP Publishing}
}

@article{hampel2022review,
  title={A review on fast tomographic imaging techniques and their potential application in industrial process control},
  author={Hampel, Uwe and Babout, Laurent and Banasiak, Robert and Schleicher, Eckhard and Soleimani, Manuchehr and Wondrak, Thomas and Vauhkonen, Marko and L{\"a}hivaara, Timo and Tan, Chao and Hoyle, Brian and others},
  journal={Sensors},
  volume={22},
  number={6},
  pages={2309},
  year={2022},
  publisher={MDPI}
}

@article{hauptmann2019real,
  title={Real-time cardiovascular MR with spatio-temporal artifact suppression using deep learning--proof of concept in congenital heart disease},
  author={Hauptmann, Andreas and Arridge, Simon and Lucka, Felix and Muthurangu, Vivek and Steeden, Jennifer A},
  journal={Magnetic resonance in medicine},
  volume={81},
  number={2},
  pages={1143--1156},
  year={2019},
  publisher={Wiley Online Library}
}

@article{kustner2020cinenet,
  title={CINENet: deep learning-based 3D cardiac CINE MRI reconstruction with multi-coil complex-valued 4D spatio-temporal convolutions},
  author={K{\"u}stner, Thomas and Fuin, Niccolo and Hammernik, Kerstin and Bustin, Aurelien and Qi, Haikun and Hajhosseiny, Reza and Masci, Pier Giorgio and Neji, Radhouene and Rueckert, Daniel and Botnar, Ren{\'e} M and others},
  journal={Scientific reports},
  volume={10},
  number={1},
  pages={13710},
  year={2020},
  publisher={Nature Publishing Group UK London}
}

@article{hakkarainen2019undersampled,
  title={Undersampled dynamic X-ray tomography with dimension reduction Kalman filter},
  author={Hakkarainen, Janne and Purisha, Zenith and Solonen, Antti and Siltanen, Samuli},
  journal={IEEE Transactions on Computational Imaging},
  volume={5},
  number={3},
  pages={492--501},
  year={2019},
  publisher={IEEE}
}

@article{tsao2003k,
  title={k-t BLAST and k-t SENSE: dynamic MRI with high frame rate exploiting spatiotemporal correlations},
  author={Tsao, Jeffrey and Boesiger, Peter and Pruessmann, Klaas P},
  journal={Magnetic Resonance in Medicine: An Official Journal of the International Society for Magnetic Resonance in Medicine},
  volume={50},
  number={5},
  pages={1031--1042},
  year={2003},
  publisher={Wiley Online Library}
}

\end{document}